\documentclass[10pt,a4paper]{article}
\usepackage[british]{babel}
\usepackage{amssymb,amsmath,amsthm}
\usepackage{cite}
\usepackage{graphicx}
\usepackage[bookmarks=true, bookmarksopen=true, bookmarksopenlevel=4]{hyperref}
\usepackage{color}
\usepackage{enumerate}
\usepackage{pgf,tikz}
\usepackage{soul}
\usetikzlibrary{decorations.pathreplacing,matrix,calc,arrows,shapes.geometric}

\hypersetup{pdftitle={Acoustic transmission problems: wavenumber-explicit bounds and resonance-free regions},
  pdfauthor={A. Moiola, E.A. Spence},
  colorlinks=true, linkcolor=myblue,  citecolor=myblue, filecolor=myblue,   urlcolor=myblue,}

\allowdisplaybreaks[4]
\definecolor{myblue}{rgb}{0,0,0.6} 
\newcommand*{\der}[2]{\frac{\partial #1}{\partial #2}}
\newcommand{\iin}{\;\text{in}\;}
\newcommand{\oon}{\;\text{on}\;}
\newcommand{\deO}{{\partial\Omega}}
\newcommand*{\conj}[1]{\overline{#1}}
\newcommand*{\N}[1]{\left\|#1\right\|}

\newcommand{\curl} {\mathop{\rm curl}\nolimits}
\newcommand{\dive} {\mathop{\rm div}\nolimits}
\def\div{\mathop{\rm div}\nolimits}

\DeclareMathOperator{\diam}{diam}
\newcommand{\uu}[1]{\hbox{\boldmath$#1$}}
\newcommand{\Uu}[1]{{\mathbf{#1}}}
\newcommand{\IN}{\mathbb{N}}\newcommand{\IR}{\mathbb{R}}\newcommand{\IZ}{\mathbb{Z}}
\newcommand{\bn}{{\Uu n}}\newcommand{\bv}{{\Uu v}}\newcommand{\bx}{{\Uu x}}\newcommand{\by}{{\Uu y}}
\newcommand{\bE}{{\Uu E}}\newcommand{\bF}{{\Uu F}}\newcommand{\bH}{{\Uu H}}
\newcommand{\bJ}{{\Uu J}}\newcommand{\bQ}{{\Uu Q}}
\newcommand{\bzero}{\Uu{0}}\newcommand{\bnu}{{\uu\nu}}  
\newtheorem{theorem}{Theorem}[section]
\newtheorem{defin}{Definition}[section]
\newtheorem{lemma}{Lemma}[section]
\newtheorem{remark}{Remark}[section]
\newtheorem{proposition}{Proposition}[section]
\newtheorem{cor}{Corollary}[section]

\newcommand{\ee}{{\rm e}}\newcommand{\ri}{{\rm i}}
\newcommand{\beq}{\begin{equation}}      \newcommand{\eeq}{\end{equation}}
\newcommand{\beqs}{\begin{equation*}}    \newcommand{\eeqs}{\end{equation*}}
\newcommand{\bit}{\begin{itemize}}       \newcommand{\eit}{\end{itemize}}
\newcommand{\ben}{\begin{enumerate}}     \newcommand{\een}{\end{enumerate}}
\newcommand{\bal}{\begin{align}}         \newcommand{\eal}{\end{align}}
\newcommand{\bals}{\begin{align*}}       \newcommand{\eals}{\end{align*}}
\newcommand{\bse}{\begin{subequations}}	 \newcommand{\ese}{\end{subequations}}
\newcommand{\bpr}{\begin{proposition}}   \newcommand{\epr}{\end{proposition}}
\newcommand{\bre}{\begin{remark}}        \newcommand{\ere}{\end{remark}}
\newcommand{\bpf}{\begin{proof}}         \newcommand{\epf}{\end{proof}}
\newcommand{\ble}{\begin{lemma}}         \newcommand{\ele}{\end{lemma}}
\newcommand{\bco}{\begin{corollary}}     \newcommand{\eco}{\end{corollary}}
\newcommand{\bex}{\begin{example}}       \newcommand{\eex}{\end{example}}
\newcommand{\bth}{\begin{theorem}}       \newcommand{\enth}{\end{theorem}}
\newcommand{\Rea}{\mathbb{R}}            \newcommand{\Com}{\mathbb{C}}
\newcommand{\Oi}{\Omega_i}

\newcommand{\pdiff}[2]{\frac{\partial #1}{\partial #2}}
\newcommand{\nus}{|u|^2}
\newcommand{\ngus}{|\nabla u|^2}
\newcommand{\nurs}{|u_r|^2}

\newcommand{\supp}{\operatorname{supp}}
\newcommand{\nvs}{|v|^2}
\newcommand{\ngvs}{|\nabla v|^2}
\newcommand{\nvrs}{|v_r|^2}
\newcommand{\gv}{\nabla v}

\newcommand{\half}{\frac{1}{2}}
\newcommand{\tendo}{\rightarrow 0}
\newcommand{\tendi}{\rightarrow \infty}
\newcommand{\dn}{{\partial_\bn}}
\newcommand{\Oin}{{\Omega_i}}
\newcommand{\Oout}{{\Omega_o}}
\newcommand{\nin}{{n_i}}
\newcommand{\nout}{{n_o}}
\newcommand{\Hankel}{H^{(1)}}
\newcommand{\gD}{{g_D}}

\newcommand{\AD}{{A_D}}
\newcommand{\AN}{{A_N}}
\newcommand{\ain}{{a_i}}
\newcommand{\aout}{{a_o}}
\newcommand{\SRC}{\mathrm{SRC}}
\newcommand{\loc}{_{\mathrm{loc}}}
\newcommand{\Rimp}{{R_{\chi}}}

\newcommand{\mythmname}[1]{\textbf{\emph{(#1.)}}}
\newcommand{\cL}{{\mathcal L}}
\newcommand{\cS}{{\mathcal S}}
\newcommand{\cO}{{\mathcal O}}
\newcommand{\cD}{{\mathcal D}}
\newcommand{\cM}{{\mathcal M}}
\newcommand{\vb}{\overline{v}}

\newcommand{\dnu}{{\partial_{\bnu}}}

\newcommand{\tin}{\text{ in }}
\newcommand{\tfa}{\text{ for all }}

\newcommand{\tas}{\text{ as }}
\newcommand{\tand}{\text{ and }}

\numberwithin{equation}{section}

\usepackage{fancyhdr}
\title{Acoustic transmission problems:\\
wavenumber-explicit bounds and resonance-free regions}

\author{Andrea Moiola\thanks{
Department of Mathematics, University of Pavia, 27100 Pavia, Italy (\texttt{andrea.moiola@unipv.it})},
\, Euan A.\ Spence\thanks{Department of Mathematical Sciences, University of Bath, Bath, BA2 7AY (\texttt{E.A.Spence@bath.ac.uk})}}

\date{\today}

\begin{document}
\maketitle

\begin{abstract}
We consider the Helmholtz transmission problem with one penetrable star-shaped Lipschitz obstacle. Under a natural assumption about the ratio of the wavenumbers, we prove bounds on the solution in terms of the data, with these bounds explicit in all parameters. 
In particular, the (weighted) $H^1$ norm of the solution is bounded by the $L^2$ norm of the source term, independently of the wavenumber.
These  bounds then imply the existence of a resonance-free strip beneath the real axis.
The main novelty is that the only comparable results currently in the literature are for smooth, convex obstacles with strictly positive curvature,
while here we assume only Lipschitz regularity and star-shapedness with respect to a point.
Furthermore, our bounds are obtained using identities first introduced by Morawetz (essentially integration by parts), whereas the existing bounds use the much-more sophisticated technology of microlocal analysis and propagation of singularities.
We also adapt existing results to show that if the assumption on the wavenumbers is lifted, then no bound with polynomial dependence on the wavenumber is possible.


\end{abstract}

\section{Introduction}\label{sec:intro}

The acoustic transmission problem, modelled by the Helmholtz equation, is a classic problem in scattering theory.
Despite having been studied from many different perspectives over the years, it remains a topic of active research. For example, recent research on this problem includes the following.
\bit
\item  
Designing novel integral-equation formulations of the transmission problem; see, e.g.,
\cite{BoDoLeTu:15} and the reviews \cite{ClHiJe:13, ClHiJePi:15}.
\item 
Designing hybrid numerical-asymptotic methods to approximate the solution of the transmission problem with a number of degrees of freedom that grows slowly (or, ideally, is constant) as the wavenumber increases \cite{GrHeLa:15,GHL17}.
\item 
Quantifying how uncertainty in the shape of the obstacle affects the solution of the transmission problem \cite{HiScScSc:18}.
\item Obtaining sharp bounds on the location of resonances of the transmission problem \cite{Ga:15}.
\item Obtaining sharp estimates on the scattered field away from the obstacle, for example in the case when the obstacle is a ball \cite{HPV07, Cap12, CLP12, NV12}, motivated by applications in imaging, inverse problems, and cloaking.
\item Designing fast solvers for the Helmholtz equation in media where the wavenumber is piecewise smooth (i.e.~transmission problems); see, e.g., the recent review \cite{GaZh:16} and the references therein. 
\item Proving wavenumber-explicit bounds on solutions of boundary value problems that approximate transmission problems, motivated by applications in numerical analysis \cite{BaChGo:16, Ch:16, SaTo:17, GrSa:18}.
\eit

In this paper we focus on the case of transmission through one obstacle; i.e.~the problem has two real wavenumbers: one inside the obstacle, and one outside the obstacle. We give a precise definition of this problem in Equation \eqref{eq:BVP} below. Our results can also be extended to more general situations (see Remark \ref{rem:extensions} below).
 
A natural question to ask about the transmission problem is:
\bit
\item[Q1.] Can one find a bound on the solution in terms of the data, with the bound explicit in the two wavenumbers?
\eit
Ideally, we would also like the bound to be either independent of the shape of the scatterer, or explicit in any of its natural geometric parameters; for example, if the domain is star-shaped, then we would ideally like the bound to be explicit in the star-shapedness parameter.
Another fundamental question is
\bit
\item[Q2.] Does the solution operator (thought of as a function of the wavenumber) have a resonance-free region underneath the real axis?
\eit
The relationship between Q1, Q2, and the question of local-energy decay for solutions of the corresponding wave equation is well-understood in scattering theory, and goes back to the work of Lax, Morawetz, and Phillips. In this particular situation of the Helmholtz transmission problem, Vodev proved in \cite[Theorem 1.1 and Lemma 2.3]{Vo:99} that an appropriate bound on the solution for real wavenumbers implies the existence of a resonance-free strip beneath the real axis.

\paragraph{Existing work on Q1 and Q2 for the transmission problem.}
To the authors' knowledge, there are five main sets of results regarding Q1 and Q2 for the Helmholtz transmission problem in the literature; we highlight that several of these results cover more general transmission problems than the single-penetrable-obstacle one considered in this paper.

\ben
\item[(a)] When the wavenumber outside the obstacle is \emph{larger} than the wavenumber inside the obstacle, and the obstacle is $C^\infty$ and convex with strictly positive curvature, Cardoso, Popov, and Vodev proved that the solution can be bounded independent of the wavenumber, and thus that there exists a resonance-free strip beneath the real axis \cite{CaPoVo:99} (these results were an improvement of the earlier work by Popov and Vodev \cite{PoVo:99a}).
\item[(b)] When the wavenumber outside the obstacle is \emph{smaller} than the wavenumber inside the obstacle, and the obstacle is $C^\infty$ and convex with strictly positive curvature, Popov and Vodev proved that there exists a sequence of complex wavenumbers (lying super-algebraically close to the real axis) through which the norm of the solution grows faster than any algebraic power of the wavenumber \cite{PoVo:99}.
\item[(c)] For either configuration of wavenumbers, and for any $C^\infty$ obstacle, Bellassoued proved that the norm of the solution cannot grow faster than exponentially with the wavenumber \cite{Be:03}.
\item[(d)] Further information about the location and the asymptotics of the resonances when the obstacle is $C^\infty$ and convex with strictly positive curvature, for both wavenumber configurations above, 
was obtained by Cardoso, Popov, and Vodev in \cite{CaPoVo:01}. 
Sharp bounds on the location of the resonances (again for both configurations of the wavenumbers) were given recently by Galkowski in \cite{Ga:15}.

\item[(e)] The case when the obstacle is a ball, in both configurations of the wavenumbers, has been studied by Capdeboscq and coauthors in \cite{CLP12,Cap12} 
(and summarised in  \cite[Chapter~5]{AC16})
using separation of variables and bounds on Bessel and Hankel functions.
\een

\paragraph{The main results of this paper and their novelty.} In this paper we prove analogues of the bound in (a) above when the obstacle is Lipschitz and star-shaped (Theorems \ref{thm:FirstBound} and \ref{thm:nonzero}) and hence also the existence of a resonance-free strip beneath the real axis (Theorem \ref{thm:res}).
Our condition on the ratio of the wavenumbers is slightly more restrictive than that in \cite{CaPoVo:99}, since the parameters in the transmission conditions are also involved (see Equation \eqref{eq:cond} and Remark \ref{rem:unnatural} below). Nevertheless
we believe this is the first time such results have been proved for the transmission problem when the obstacle is non-convex or non-smooth.

Furthermore, the constant in our bound is given completely explicitly, and the bound is valid for all wavenumbers greater than zero (satisfying the restriction on the ratio).
On the other hand, the bound for smooth convex obstacles in \cite{CaPoVo:99} assumes that both wavenumbers are large, with the ratio of the two fixed, but the bound is not explicit in this ratio (although we expect the results of \cite{Ga:15} could be used to get a bound for smooth convex obstacles that is explicit in the ratio of wavenumbers, when the wavenumbers are large enough \cite{Ga:16}).
An additional feature of the constant in our bound is that, whilst the bound is valid for all star-shaped Lipschitz obstacles, the constant only depends on the diameter of the obstacle. This ``shape-robustness" makes the bound particular suitable for applications in quantifying how uncertainty in the shape of the obstacle affects the solution (done for the transmission problem at low frequency in \cite{HiScScSc:18}), and these applications are currently under investigation.

We highlight that the bound in \cite{CaPoVo:99} relies on microlocal analysis and the deep results of Melrose and Sj\"ostrand on propagation of singularities \cite{MeSj:78, MeSj:82}. 
In contrast, our bound is obtained using identities for solutions of the Helmholtz equation first introduced by Morawetz in \cite{MoLu:68, Mo:75}, which boil down to multiplying the PDEs by carefully-chosen test functions and integrating by parts.
Whereas Morawetz's identities have been used to prove bounds on many Helmholtz BVPs, famously the exterior Dirichlet and Neumann problems in \cite{MoLu:68, Mo:75}, it appears that (surprisingly) this paper is their first application to the Helmholtz transmission problem involving one penetrable obstacle (Remark \ref{rem:rough} below discusses their applications to transmission problems not involving a bounded obstacle). 
The novelty of this paper is therefore not in the techniques that are used, but the fact that these well-known techniques can be applied to a classic problem to obtain new results.

\paragraph{Outline of the paper.}
In \S\ref{sec:form} we define the Helmholtz transmission problem for Lipschitz obstacles and recap results on existence, uniqueness, and regularity.
In \S\ref{sec:result} we give the main results (Theorems~\ref{thm:FirstBound}, \ref{thm:nonzero} and \ref{thm:res}).
In \S\ref{sec:Morawetz} we derive the Morawetz identities used in the proofs of Theorems~\ref{thm:FirstBound} and \ref{thm:nonzero}, and in \S\ref{sec:proofs} we prove the main results.
In \S\ref{sec:blowup} we adapt the existing results of \cite{PoVo:99} about super-algebraic growth of the norm of the solution through a sequence of complex wavenumbers to prove analogous growth through a sequence of real wavenumbers; we illustrate this growth through real wavenumbers with plots when the obstacle is a 2-d ball.

\paragraph{Motivation for \S\ref{sec:blowup}.}
Our motivation for adapting the results of \cite{PoVo:99} (and also highlighting the results of \cite{Cap12, CLP12}) in \S\ref{sec:blowup} is the recent investigations 
\cite{Ch:16, BaChGo:16, SaTo:17, GrSa:18} of the interior impedance problem for the Helmholtz equation with piecewise-constant wavenumber (and the related investigation  \cite{OhVe:16} for piecewise-Lipschitz wavenumber).
These investigations, coming from the numerical-analysis community, concern Helmholtz transmission problems, but consider the interior impedance problem, because the impedance boundary condition is a simple approximation of the Sommerfeld radiation condition that is commonly used when implementing finite-element methods (see Remark~\ref{rem:trunc} below).

The results of \cite{PoVo:99}, adapted as in \S\ref{sec:blowup} (using, in particular, recent wavenumber-explicit bounds on layer-potential operators from \cite{Sp2013a, HaTa:15}),
\bit
\item[(a)] show that the ``technical'' assumptions on the wavenumber in \cite[\S1]{BaChGo:16} are in fact necessary for the results of  \cite[\S1]{BaChGo:16} to hold, and 
\item[(b)] partially answer the conjecture in \cite[\S2.3]{SaTo:17} about the maximal growth of the norm of the solution operator. 
\eit

\section{Formulation of the problem}\label{sec:form}
\subsection{Geometric notation.}\label{sec:Notation}

Let $\Oin\subset\IR^d$, $d \geq 2$, be a bounded Lipschitz open set.
Denote $\Oout:=\IR^d\setminus\overline\Oin$ and $\Gamma:=\deO_i=\deO_o$.
Let $\bn$ be the unit normal vector field on $\Gamma$ pointing from $\Oin$ into $\Oout$.
We denote by $\dn$ the corresponding Neumann trace from one of the two domains $\Oin$ and $\Oout$ and we do not use any symbol for the Dirichlet trace on $\Gamma$.
For any $\varphi\in L^2\loc(\IR^d)$, we write $\varphi_i:=\varphi|_\Oin$ and $\varphi_o:=\varphi|_\Oout$.

For $a>0$ and $\bx_0\in\IR^d$ we denote by $B_a(\bx_0)=\{\bx\in\IR^d:|\bx|<a\}$ the ball with centre $\bx_0$ and radius $a$; if $\bx_0=\bzero$ we write $B_a=B_a(\bzero)$.
Given $R>0$ such that $\conj\Oin\subset B_R$,
let  $D_R:=B_R\cap\Oout$ and $\Gamma_R:=\partial B_R=\{|\bx|=R\}$.
On $\Gamma_R$ the unit normal $\bn$ points outwards.
With $D$ denoting an open set or a $d-1$-dimensional manifold, $\N{\cdot}_{D}$ denotes $L^2(D)$ norm for scalar or vector fields.
On $\Gamma$ and $\Gamma_R$, $\nabla_T$ denotes the tangential gradient.

To state the main results, we need to define the notions of  \emph{star-shaped} and \emph{star-shaped with respect to a ball}. 

\begin{defin}
(i) $\Oin$ is \emph{star-shaped with respect to the point $\bx_0$} if, whenever $\bx \in \Oin$, the segment $[\bx_0,\bx]\subset \Oin$.

\noindent (ii) $\Oin$ is \emph{star-shaped with respect to the ball $B_{a}(\bx_0)$} if it is star-shaped with respect to every point in $B_{a}(\bx_0)$.
\end{defin}

These definitions make sense even for non-Lipschitz $\Oin$, but when $\Oin$ is Lipschitz one can characterise star-shapedness with respect to a point or ball in terms of $(\bx-\bx_0)\cdot \bn(\bx)$ for $\bx\in\Gamma$.
\begin{lemma}\label{lem:star}
(i) If $\Oin$ is Lipschitz, then it is star-shaped with respect to $\bx_0$ if and only if $(\bx-\bx_0)\cdot\bn(\bx)\geq 0$ for all $\bx \in\Gamma$ for which $\bn(\bx)$ is defined.

\noindent (ii) $\Oin$ is star-shaped with respect to $B_{a}(\bx_0)$ if and only if 
it is Lipschitz and
$(\bx-\bx_0) \cdot \bn(\bx) \geq {a}$ for all  $\bx \in \Gamma$ for which $\bn(\bx)$ is defined; 
\end{lemma}

\bpf 
See \cite[Lemma 5.4.1]{AndreaPhD} or \cite[Lemma 3.1]{MaxwellPDE}.
\epf

In the rest of the paper, whenever $\Oin$ is star-shaped with respect to a point or ball, we assume (without loss of generality) that $\bx_0=\bzero$.

\subsection{The Helmholtz transmission problem.}

From the point of view of obtaining wavenumber-explicit bounds, we are interested in the case when the wavenumber is real. Nevertheless, in order to talk about resonance-free regions, we must also consider complex wavenumbers.

\begin{defin}\mythmname{Sommerfeld radiation condition}
Given $\varphi\in C^1(\IR^d\setminus B_R)$, for some ball $B_R=\{|\bx|<R\}$,  and $\kappa\in \Com\setminus\{0\}$ with $\Im \kappa\geq 0$, we say that $\varphi$ satisfies the Sommerfeld radiation condition if 
\beq
\label{eq:src}
\lim_{r\to\infty}r^{\frac{d-1}2}\left(\der{\varphi(\bx)}r-\ri\kappa \varphi(\bx)\right)=0
\eeq
uniformly in all directions, where $r=|\bx|$; we then write $\varphi\in\SRC(\kappa)$.
\end{defin}
Recall that when $\Im \kappa>0$, if $\varphi$ satisfies the Sommerfeld radiation condition, then $\varphi$ decays exponentially at infinity (see, e.g., \cite[Theorem 3.6]{CoKr:83}).

\begin{defin}\mythmname{The Helmholtz transmission problem}\label{def:HTP}
Let $k\in \Com\setminus\{0\}$ with $\Im k\geq 0$, and let $\nin,\nout,\ain,\aout,\AD,\AN$ be positive real numbers.
Let $f_i\in L^2(\Oin)$, $f_o\in L^2(\Oout)$, $g_D\in H^1(\Gamma)$, $g_N\in L^2(\Gamma)$, and assume $f_o$ has compact support.
The Helmholtz transmission problem is: find $u\in H^1\loc(\IR^d\setminus\Gamma)$ such that,
\begin{align}
\begin{aligned}
\ain\Delta u_i+k^2 \nin u_i &=f_i &&\iin\Oin,\\
\aout\Delta u_o+k^2 \nout u_o &=f_o &&\iin\Oout,\\
u_o&=\AD u_i+g_D &&\oon \Gamma,\\
\aout\dn u_o&=\AN\ain\dn u_i+g_N &&\oon \Gamma,\\ 
&&&u_o\in\SRC(k\sqrt{n_o/a_o}).
\end{aligned}
\label{eq:BVP}
\end{align}
\end{defin}
Four of the parameters $\nin,\nout,\ain,\aout,\AD,\AN$ are redundant; in particular we can set either $\AD=\ain=\aout=\nout=1$ or $\AD=\AN=\aout=\nout=1$ and still cover all problems by rescaling the remaining coefficients, $u_i$, and the source terms.
Nevertheless, we keep all six parameters in \eqref{eq:BVP} since given a specific problem it is then easy to write it in the form \eqref{eq:BVP}, setting some parameters to one.

Some extensions to transmission problems more general than those in Definition~\ref{def:HTP} are discussed in Remark~\ref{rem:extensions}.

\bre\mythmname{Relation to acoustics and electromagnetics}\label{rem:physics}
Time-harmonic acoustic transmission problems are often written in the form
\begin{align}
\div\Big(\frac1\rho\nabla u\Big)+\frac{\kappa^2}\rho u=F,
\qquad
u\in H^1(\IR^d),\qquad \frac1\rho\nabla u\in H(\dive;\IR^d),\quad
\label{eq:BVPrho}
\end{align}
and $u$ satisfies the Sommerfeld radiation condition, where $\rho(\bx)$ and $\kappa(\bx)$ are positive functions; see e.g.\ \cite[eq.~(1)]{HMK02}.
(Recall that $\bv\in H(\dive;\IR^d)$ if and only if $\bv|_\Oin\in H(\dive;\Oi)$, $\bv|_\Oout\in H(\dive;\Oout)$ and $\bv|_\Oin\cdot\bn=\bv|_\Oout\cdot\bn$ in  $H^{-1/2}(\Gamma)$.)
In the particular case where $\rho$ and $\kappa$ take two different values on $\Oin$ and $\Oout$, problem \eqref{eq:BVPrho} can be written in the form \eqref{eq:BVP} choosing, for example, 
\[
\AD=\AN=1,\quad
a=\frac1\rho,
\quad  k = \kappa_o, \quad  \nout=\frac1\rho_o,\quad
\nin=\Big(\frac{\kappa_i}{\kappa_o}\Big)^2\frac1\rho_i,\quad f=F,\quad g_D=g_N=0,
\]
or
\[
\AD=\ain=\aout=\nout=1,\quad
\AN=\frac{\rho_o}{\rho_i},
\quad  k = \kappa_o, \quad  
\nin=\Big(\frac{\kappa_i}{\kappa_o}\Big)^2,\quad f=\rho F,\quad g_D=g_N=0.
\]
(More generally, one can choose any constant $k>0$ and then let $n=a\kappa^2/ k^2$.)

The time-harmonic Maxwell equations are
\begin{align}\label{eq:Maxwell}
\curl \bH+\ri k\varepsilon\bE=(\ri k)^{-1}\bJ,\qquad
\curl\bE-\ri k\mu\bH=\bzero
\qquad \iin \IR^3.
\end{align}
When all fields and parameters involved depend only on two Cartesian space variables, say $x$ and $y$, Equations \eqref{eq:Maxwell}
 reduce to the (heterogeneous) Helmholtz equation in $\IR^2$.
In the transverse-magnetic (TM) mode, $\bJ$ and $\bE$ are given by $\bJ=(0,0,J_z(x,y))$ and $\bE=(0,0,$ $E_z(x,y))$, so \eqref{eq:Maxwell} reduce to a scalar equation for the third component of the electric field:
\begin{align*}
\div\left(\frac1\mu\nabla E_z\right)+k^2\varepsilon E_z=-J_z.
\end{align*}
If the permittivity $\varepsilon$ and the permeability $\mu$ are constant in $\widetilde\Oin=\Oin\times\IR\subset\IR^3$ and $\widetilde\Oout=\Oout\times\IR\subset\IR^3$
for $\Oin,\Oout\subset\IR^2$ as in \S\ref{sec:Notation}, then \eqref{eq:Maxwell} (supplemented with suitable radiation conditions) can be written as Problem \eqref{eq:BVP}
for $u=E_z$ in $\IR^2$ with $a=1/\mu$, $n=\varepsilon$, $\AD=\AN=1$.
Vice versa, in the transverse-electric (TE) mode, $\bJ=(J_x(x,y),J_y(x,y),0)$, $\bH=(0,0,H_z(x,y))$ and 
\begin{align*}
\div\left(\frac1\varepsilon\nabla H_z\right)+k^2\mu H_z=(\ri k\varepsilon)^{-1}\Big(\der{J_x}y-\der{J_y}x\Big),
\end{align*}
so \eqref{eq:Maxwell} can be written as \eqref{eq:BVP} for $u=H_z$ with $a=1/\varepsilon$, $n=\mu$, $\AD=\AN=1$.
Observe that for TM and TE modes, the parameters $\nin, \nout, \ain, \aout$ depend on properties of the medium through which the waves propagate, whereas $k$ depends on the wave itself.
\ere

\begin{defin}\mythmname{Scattering problem}
Let $k\in \Com\setminus\{0\}$ with $\Im k\geq 0$, and let $\nin,\nout,\ain$, $\aout$, $\AD$, $\AN$ be positive real numbers.
Let $u^I$ be a solution of $\aout\Delta u^I+k^2\nout u^I=0$ that is $C^\infty$ in a neighbourhood of $\conj\Oin$ (for example a plane wave, a circular or spherical wave, or a fundamental solution centred in $\Oout$).
Define the total field $u^T$ to be solution of 
\begin{align}
\begin{aligned}
\ain\Delta u_i^T+k^2 \nin u_i^T &=0 &&\iin\Oin,\\
\aout\Delta u_o^T+k^2 \nout u_o^T &=0 &&\iin\Oout,\\
u_o^T&=\AD u_i^T &&\oon \Gamma,\\
\aout\dn u_o^T&=\AN\ain\dn u_i^T &&\oon \Gamma,\\ 
&&& (u_o^T-u^I)\in\SRC(k\sqrt{\nout/\aout}).
\end{aligned}
\label{eq:BVPuT}
\end{align}
The scattered field defined by $u:=u^T-u^I$ satisfies
\begin{align}
\begin{aligned}
\ain\Delta u_i+k^2 \nin u_i &=k^2\Big(\frac\ain\aout\nout-\nin\Big)u^I=:f_i &&\iin\Oin,\\
\aout\Delta u_o+k^2 \nout u_o &=0 &&\iin\Oout,\\
u_o&=\AD u_i +(\AD-1)u^I&&\oon \Gamma,\\
\aout\dn u_o&=\AN\ain\dn u_i+(\AN-1)\ain\dn u^I &&\oon \Gamma,\\ 
&&& u_o\in\SRC(k\sqrt{\nout/\aout}).
\end{aligned}
\label{eq:BVPuS}
\end{align}
The scattering problem \eqref{eq:BVPuS} can therefore be written in the form \eqref{eq:BVP} for $f_i=k^2(\frac\ain\aout\nout-\nin)u^I$, $f_o=0$, $g_D=(\AD-1)u^I$ and $g_N=(\AN-1)\ain\dn u^I $.
\end{defin}

The next lemma, proved in Appendix \ref{sec:appA}, addresses the questions of existence, uniqueness, and regularity of the solution of \eqref{eq:BVP}.

\begin{lemma}\mythmname{Existence, uniqueness, and regularity}\label{lem:exist}
The Helmholtz transmission problem of Definition \ref{def:HTP} admits a unique solution $u\in H^1\loc(\IR^d\setminus\Gamma)$.
Moreover $u_i,u_o\in H^1(\Gamma)$ and $\dn u_i,\dn u_o\in L^2(\Gamma)$.
\end{lemma}

The key point about Lemma \ref{lem:exist} is that, whilst the existence and uniqueness results are well known, the 
regularity results $u_i,u_o\in H^1(\Gamma)$ and $\dn u_i,\dn u_o\in L^2(\Gamma)$ are currently only available in the literature in the case $d=3$, $f_i = f_o=0$, and $k\in \Rea$. 
These results are consequences of the harmonic-analysis results about layer potentials in \cite{CoMcMe:82, Ve:84, EsFaVe:92}
and the regularity results of Ne\v{c}as for strongly elliptic systems in \cite[\S5.1.2 and \S5.2.1]{Ne:67}, \cite[Theorem 4.24]{MCL00}, and are needed to apply the Morawetz identities of \S\ref{sec:Morawetz} to the solution of the transmission problem when $\Oi$ is Lipschitz.

Finally, recall that $k$ is a resonance of the boundary value problem \eqref{eq:BVP} if there exists a non-zero $u_o$ satisfying \eqref{eq:BVP} with $f_i=f_o=g_D=g_N=0$ and the Sommerfeld radiation condition replaced by 
\beq\label{eq:BU1}
u_o (\bx) = \frac{{\ee}^{\ri k\sqrt{n_o/a_o}r}}{r^{(d-1)/2}}\left( u_\infty(\widehat{\bx}) + \cO\left(\frac{1}{r}\right)\right) \quad\tas r:= |\bx|\tendi,
\eeq
for some function $u_\infty$ of $\widehat{\bx}=\bx/|\bx|$ (the far-field pattern); see, e.g., \cite[\S3.6 and Theorem 4.9]{DyZw:16} or \cite[\S2]{PaK09}.
The uniqueness result of Lemma \ref{lem:exist} implies that any resonance must have $\Im k<0$, and thus \eqref{eq:BU1} implies that $u_o$ grows exponentially at infinity.

\section{Main results}\label{sec:result}

\subsection{Bounds on the solution (answering Q1)}
\label{sec:resultQ1}

In this section we assume that $k>0$, but analogous results hold for $k<0$ as well.
We recall that the main assumptions we stipulate, namely \eqref{eq:cond} and \eqref{eq:gDgNcond} below, mean that the wavenumber is larger in the exterior region $\Oout$ than in the interior region $\Oin$, or equivalently that the wavelength $\lambda=\frac{2\pi\sqrt{a}}{k\sqrt n}$ is longer in $\Oin$ than in $\Oout$.
We recall from \S\ref{sec:Notation} that the notation $\N{\cdot}_D$ stands for the $L^2(D)$ norm.

\begin{theorem}\label{thm:FirstBound}
Assume that $\Oi$ is star-shaped,
\beq\label{eq:cond}
\frac{\nin}{\nout}\leq \frac{A_D}{A_N}\leq \frac{a_i}{a_o},
\eeq
$g_N=g_D=0$ and $k>0$. 
Given $R>0$ such that $\supp f_o \subset B_R$, recall that $D_R:= \Oout \cap B_R$. The solution of BVP \eqref{eq:BVP} then satisfies
\begin{align}
\begin{aligned}\label{eq:bound1}
&a_i\N{\nabla u_i}_{\Oin}^2+k^2n_i\N{u_i}^2_{\Oin}
+\frac1{A_D A_N}\left(a_o\N{\nabla u_o}_{D_R}^2+k^2n_o\N{u_o}^2_{D_R}\right)\\
&\hspace{2cm}\le
\bigg[\frac{4\diam(\Omega_i)^2}{a_i}+\frac{1}{n_i}\left(2\sqrt{\frac{n_o}{a_o}} R +\frac{d-1}{k}\right)^2\bigg]
\N{f_i}_{\Omega_i}^2\\
&\hspace{30mm}+\frac1{A_D A_N}
\bigg[\frac{4R^2}{a_o}+\frac{1}{\nout}\left(2\sqrt{\frac{n_o}{a_o}} R +\frac{d-1}{k}\right)^2\bigg]\N{f_o}_{D_R}^2.
\end{aligned}
\end{align}
\end{theorem}
The bound \eqref{eq:bound1} is valid for all star-shaped Lipschitz $\Omega_i$, but the constants on the right-hand side only depend on $\Oi$ via $\diam(\Omega_i)$. 
As highlighted in \S\ref{sec:intro}, we expect that this uniformity of the bound with respect to the geometry makes it particular suitable for applications in quantifying how uncertainty in the shape of $\Oi$ affects the solution (as done for small $k$ in \cite{HiScScSc:18}).

In Theorem \ref{thm:FirstBound} we assumed that the boundary source terms $g_D$ and $g_N$ vanish.
In the next theorem we consider general $g_D\in H^1(\Gamma)$ and $g_N\in L^2(\Gamma)$.
In order to do this, we need to assume that the inequalities \eqref{eq:cond} on the parameters are strict and that $\Oin$ is star-shaped with respect to a ball.

\begin{theorem}\label{thm:nonzero}
Assume that $\Oin$ is star-shaped with respect to $B_{\gamma\diam(\Oin)}$ for some $0<\gamma\le1/2$, 
\begin{equation}\label{eq:gDgNcond}
\frac\nin\nout<\frac\AD\AN<\frac\ain\aout,
\end{equation}
$k>0$ and $R>0$ is such that $\supp f_o\subset B_R$.
Then the solution of \eqref{eq:BVP} satisfies
\begin{align}\label{eq:gDgNbound}
\begin{aligned}
&a_i\N{\nabla u_i}_{\Oin}^2+k^2n_i\N{u_i}^2_{\Oin}
+\frac1{A_D A_N}\left(a_o\N{\nabla u_o}_{D_R}^2+k^2n_o\N{u_o}^2_{D_R}\right)\\
&\hspace{0mm}\le
\bigg[\frac{4\diam(\Omega_i)^2}{a_i}+\frac{1}{n_i}\left(2\sqrt{\frac{n_o}{a_o}} R +\frac{d-1}{k}\right)^2\bigg]
\N{f_i}_{\Omega_i}^2\\
&\hspace{7mm}+\frac1{A_D A_N}
\bigg[\frac{4R^2}{a_o}+\frac{1}{\nout}\left(2\sqrt{\frac{n_o}{a_o}} R +\frac{d-1}{k}\right)^2\bigg]\N{f_o}_{D_R}^2\\
&\hspace{7mm}+2\bigg[\frac{\diam(\Oin)\aout\big((3+2\gamma) \ain \AN + 2 \aout\AD\big)}{\AD\AN\gamma(\ain\AN-\aout\AD)}\bigg]\N{\nabla_T\gD}_\Gamma^2\\
&\hspace{7mm}
+2\bigg[\frac{2\diam(\Oin) n_o^2}{\gamma\AN(\nout\AD-\nin\AN)}
+\frac{(3+\gamma)\ain\Big(\nout R^2+\frac{\aout(d-1)^2}{4k^2}\Big)}{\gamma\AD\diam(\Oin)(\ain\AN-\aout\AD)}\bigg]
k^2\N{\gD}_\Gamma^2
\\
&\hspace{7mm}+\frac{2}{\gamma\aout\AN\AD}\bigg[
\frac{\diam(\Oin)(4\ain\AN+2\aout\AD)}{\ain\AN-\aout\AD}
+\frac{2\AD\Big(\nout R^2+\frac{\aout(d-1)^2}{4k^2}\Big)}{\diam(\Oin)(\nout\AD-\nin\AN)}
\bigg]\N{g_N}_\Gamma^2.
\end{aligned}
\end{align}
\end{theorem}
Note that each of the coefficients in front of the norms on the right-hand side of the bound \eqref{eq:gDgNbound} is a 
non-increasing function of $k$, apart from the coefficient $k^2$ multiplying $\N{\gD}_\Gamma^2$.

\begin{cor}\mythmname{Theorems~\ref{thm:FirstBound} and \ref{thm:nonzero} applied to the scattering problem \eqref{eq:BVPuS}}
\label{cor:ScatteringBound}
The solution $u$ of the scattering problem \eqref{eq:BVPuS} with $\nout=\aout=\AD=\AN=1$, $\nin\le1\le\ain$, and $\Oi$ star-shaped satisfies
\begin{align*}
\ain\N{\nabla u_i}_{\Oin}^2+k^2n_i\N{u_i}^2_{\Oin}
&+\N{\nabla u_o}_{D_R}^2+k^2\N{u_o}^2_{D_R}\\
&\le\bigg[\frac{4\diam(\Omega_i)^2}\ain+\frac{1}{n_i}\left(2 R +\frac{d-1}{k}\right)^2\bigg]
k^4(\ain-\nin)^2\N{u^I}_{\Omega_i}^2.
\end{align*}
The solution $u$ of the scattering problem \eqref{eq:BVPuS} with $\nout=\aout=\ain=\AD=1$, $\nin<1/\AN<1$ and $\Oi$ star-shaped with respect to $B_{\gamma\diam(\Oin)}$ (for some $0<\gamma\le1/2$) satisfies
\begin{align*}
\N{\nabla u_i}_{\Oin}^2+&k^2n_i\N{u_i}^2_{\Oin}
+\frac1\AN\Big(\N{\nabla u_o}_{D_R}^2+k^2\N{u_o}^2_{D_R}\Big)\\
&\le\bigg[4\diam(\Omega_i)^2+\frac{1}{n_i}\left(2 R +\frac{d-1}{k}\right)^2\bigg]
k^4(1-\nin)^2\N{u^I}_{\Omega_i}^2\\
&\quad+\frac4{\gamma\AN}\bigg[
\frac{\diam(\Oin)(2\AN+1)}{\AN-1}
+\frac{R^2+\frac{(d-1)^2}{4k^2}}{\diam(\Oin)(1-\nin\AN)}
\bigg](\AN-1)^2\N{\dn u^I}_\Gamma^2.
\end{align*}
\end{cor}

\bre\mythmname{Extensions of Theorems \ref{thm:FirstBound} and \ref{thm:nonzero}}\label{rem:extensions}
We have only considered the case of a single penetrable obstacle with $a_i, a_o, n_i$, and $n_o$ all real and constant, but analogues of Theorems~\ref{thm:FirstBound} and \ref{thm:nonzero}  hold in the following cases (and also when the cases are combined).
\ben
\item When there are multiple ``layers'', each with constant $a$ and $n$, and with the boundaries of the layers star-shaped with respect to the origin (for the analogue of Theorem \ref{thm:FirstBound}) or star-shaped with respect to balls centred at the origin (for the analogue of Theorem \ref{thm:nonzero}); in this case the conditions \eqref{eq:cond}/\eqref{eq:gDgNcond} must hold at each interface.
\item When $a_i, a_o, n_i$, and $n_o$ are functions of position and satisfy conditions that ensure nontrapping of rays.
\item When $\Oin$ contains an impenetrable star-shaped Dirichlet scatterer.
\item When $\Oout$ is truncated by a star-shaped boundary and the radiation condition is approximated by an impedance boundary condition.
\item When $n_i$ is complex with $0<\Im n_i\leq \delta/k$, for $\delta$ a sufficiently small constant; this models a particular case of a lossy scatterer in a lossless background.
\item When Condition \eqref{eq:cond} is partly violated, namely when $\nin/\nout$ is slightly larger (in a $k$-dependent way) than $\AD/\AN$ and $\ain/\aout$.
\een
The extension to Case 1 is clear from the proofs in \S\ref{sec:proofs}, the extension to Case 2 and 3 are covered in  \cite{GrPeSp:18},  the extension to the Case 4 is discussed in Remark~\ref{rem:impbc}, the extension to Case 5 is discussed in Remarks \ref{rem:complex} and \ref{rem:complex2}, and the extension to Case 6 is described in Proposition~\ref{thm:BoundViaTrace2}.
\ere

\subsection{Resonance-free strip (answering Q2)}

We let $R(k)$ denote the solution operator of the Helmholtz transmission problem of Definition \ref{def:HTP} when $g_D$ and $g_N$ are both zero, i.e.
\beqs
R(k): \left(
\begin{array}{c}
f_i\\
f_o
\end{array}
\right)
\mapsto
\left(
\begin{array}{c}
u_i\\
u_o
\end{array}
\right).
\eeqs
Although $R(k)$ depends also on the parameters $n_o, n_i, a_o, a_i, A_D,$ and $A_N$, in what follows we consider these fixed and consider $k$ as variable.
Let $\chi_1, \chi_2 \in C_0^\infty(\Rea^d)$ such that $\chi_j\equiv 1$ in a neighbourhood of $\Oi$, and let 
\beq\label{eq:Rimp}
\Rimp(k):= \chi_1 R(k) \chi_2;
\eeq
i.e.~$\Rimp(k)$ is the cut-off resolvent. Then
\beqs
\Rimp(k): L^2(\Oin)\oplus L^2(\Oout) \rightarrow H^1(\Oin)\oplus H^1(\Oout)
\eeqs
for $k\in \Rea\setminus\{0\}$.

\bth\mythmname{Pole-free strip beneath the real axis}\label{thm:res} 
The operator family $\Rimp(k)$ defined above is holomorphic on $\Im k>0$. 
Assume that $\Oi$ is star-shaped and the condition \eqref{eq:cond} is satisfied.
Then, there exists $C_j>0$, $j=1,2,3$ such that $\Rimp(k)$
extends from the upper-half plane to a holomorphic operator family on 
$|\Re k|\geq C_1  , \Im k\geq -C_2$
satisfying the  estimate 
\beq\label{eq:55}
\N{\Rimp(k)}_{L^2(\Oin)\oplus L^2(\Oout) \rightarrow L^2(\Oin)\oplus L^2(\Oout)} \leq\frac{C_3}{|k|}
\eeq 
in this region.
\end{theorem}

This follows from the bound of Theorem \ref{thm:FirstBound} using the result \cite[Lemma 2.3]{Vo:99}. 
Recall that this result of Vodev takes a resolvent estimate on the real axis and converts it into a resolvent estimate in a strip beneath the real axis.
In principle, one could go into the details of this result and make the width of the strip explicit in the constant from the bound on the real axis. Since  Theorem \ref{thm:FirstBound}  gives an explicit expression for that constant, we would then have an explicit lower bound for the width of the strip.

\subsection{Discussion of the main results in the context of previous results}\label{sec:previous}

We now discuss Theorems \ref{thm:FirstBound}--\ref{thm:nonzero} and \ref{thm:res} in the context of the results summarised in (a)--(e) of \S\ref{sec:intro}
(i.e.~\cite{PoVo:99,PoVo:99a,CaPoVo:99,CaPoVo:01,Be:03,Ga:15,CLP12,Cap12}) and other related work. 

We focus on results about the Helmholtz transmission problem of Definition \ref{def:HTP}, i.e.~one penetrable obstacle and piecewise-constant wavenumber. Many of these results apply to the more-general case when the wavenumber is piecewise-smooth, but we focus on the piecewise-constant case. 
There is also a substantial literature on the Helmholtz equation with continuous wavenumber, including \cite{BlKa:77} and \cite{PeVe:99}; for a survey of these result we refer the reader to \cite[\S2.4]{GrPeSp:18}. 
At the end of this section we briefly discuss (i) truncated Helmholtz transmission problems (in Remark~\ref{rem:trunc}), (ii) Helmholtz transmission problems with piecewise-constant wavenumber but \emph{not} involving a bounded obstacle (in Remark \ref{rem:rough}), and (iii) Helmholtz transmission problems when $\nin\in \Com$ with $\Im \nin >0$ (in Remark \ref{rem:complex}).

The results summarised in (a)--(d) of \S\ref{sec:intro} all consider the case when
$a_i=a_o=n_o=A_D=1$ and $g_D=g_N=0$; that is, the BVP
\begin{align}
\begin{aligned}
\Delta u_i+k^2 n_i u_i &=f_i &&\iin\Oin,\\
\Delta u_o+k^2 u_o &=f_o &&\iin\Oout,\\
u_o&=u_i &&\oon \Gamma,\\
\dn u_o&=\AN\dn u_i &&\oon \Gamma,\\ 
&&&u_o\in\SRC(k).
\end{aligned}\label{eq:BVP2}
\end{align}
The results summarised in (e) of \S\ref{sec:intro} consider \eqref{eq:BVP2} with $A_N=1$; these results are in a slightly different direction to those of (a)--(d) (involving bounds in different norms) and so we discuss them separately in Remark \ref{rem:Capdeboscq} below.

In this discussion we use the notation for the cut-off resolvent $\Rimp(k)$, and observe that the bound \eqref{eq:bound1} in the case of the BVP \eqref{eq:BVP2} is essentially equivalent to the following bound: given $k_0>0$ there exist $C_m, m=0,1,$ such that
\beq\label{eq:res1}
\N{\Rimp(k)}_{L^2 \rightarrow H^m} \leq C_m k^{m-1} \quad\tfa k\geq k_0,
\eeq
where $L^2$ and $H^0$ denote $L^2(\Oin)\oplus L^2(\Oout)$, $H^1$ denotes $H^1(\Oin)\oplus H^1(\Oout)$, and the constants $C_m$ are given explicitly in terms of $n_i, A_N$, $d$, $\diam(\Oin)$, and $R$ (with $R$ such that the support of the cut-off function $\chi_1$ appearing in the definition of $\Rimp(k)$ \eqref{eq:Rimp} is contained in $B_R$).
Once the bound 
\beq\label{eq:res2}
\N{\Rimp(k)}_{L^2 \rightarrow L^2} \leq C_0 k^{-1} \quad\tfa k\geq k_0,
\eeq
is proven, Green's identity can be used to prove the $L^2 \rightarrow H^1$ bound, with the constant $C_1$ in \eqref{eq:res1} then given explicitly in terms of 
$n_i, A_N$, $d$, $\diam(\Oin)$, $R$, and $C_0$ (the analogous argument for scattering by impenetrable Dirichlet or Neumann obstacles is given in, e.g., \cite[Lemma~2.2]{Sp2013a}).

Cardoso, Popov, and Vodev \cite{CaPoVo:99} proved the bound \eqref{eq:res2} 
when $\Oi$ is a smooth, convex obstacle with strictly positive curvature; the existence of a resonance-free strip then followed from Vodev's result in \cite{Vo:99}.
The bound in \cite{CaPoVo:99} is proved under the assumptions $n_i<1$ and $A_N>0$ and
the dependence of the constant $C_0$ on $n_i$ and $A_N$ is not given. 
The conditions $n_i<1$ and $A_N>0$ are less restrictive than our condition \eqref{eq:cond}, which in this situation is $n_i \leq 1/A_N\le 1$ 
(see Remark \ref{rem:unnatural} below for how this condition appears in our proof).

The particular case when $\Oi$ is a ball shows that a strip is the largest region one can prove is free of resonances in the case $n_i<1$ and $A_N>0$.
This was known in \cite{PoVo:99a}, but the recent results of Galkowski \cite{Ga:15} in the case when $\Oi$ is $C^\infty$ with strictly positive curvature include bounds on the width of the resonance-free strip in terms of appropriate averages of the reflectivity and chord lengths of the billiard ball trajectories in $\Oi$ \cite[Theorem 1]{Ga:15}, and these bounds are sharp when $\Oi$ is a ball \cite[\S12]{Ga:15}. 
Furthermore, these results (which build on earlier work by Cardoso, Popov, and Vodev\cite{CaPoVo:01}) show that if $\nin <1$ and $\sqrt{\nin}<1/A_N$, then the resonances themselves lie in a strip (i.e.~there exist $C_1, C_2>0$ such that the resonances $k$ satisfy $-C_2 \leq \Im k \leq -C_1$).

Popov and Vodev \cite{PoVo:99} showed that when $n_i>1$ and $A_N>0$ there exists a sequence of resonances tending to the real-axis and one has super-algebraic growth of $\|\Rimp(k)\|_{L^2\rightarrow L^2}$ through a sequence of complex wavenumbers with super-algebraically small imaginary parts; we recap this result in more detail in \S\ref{sec:PoVo}. Exponential growth in $k$ is the fastest growth possible by the results of Bellassoued \cite{Be:03}. Indeed, he proved that for any $n_i>0$ and any $A_N>0$ there exist $\widetilde{C_j}$, $j=1,\ldots,5$, such that
\beqs
\N{\Rimp(k)}_{L^2\rightarrow L^2} \leq \widetilde{C}_1 \exp(\widetilde{C}_2 |\Re k|)
\eeqs
in the region
\beqs
\Re k\geq \widetilde{C}_3, \quad \Im k \geq - \widetilde{C}_4 \exp(-\widetilde{C}_5 |\Re k|),
\eeqs
implying there is always an exponentially-small region free of resonances.

\bre\mythmname{Transmission problems when the obstacle is a  2- or 3-d ball}\label{rem:Capdeboscq}
When $\Oin$ is a ball, the solution of BVP \eqref{eq:BVP} can be written explicitly using separation of variables and expansions in Fourier--Bessel functions.
Capdeboscq and co-authors considered this problem for BVP \eqref{eq:BVP2} with $\AN=1$ for $d=2$ in \cite{Cap12} and for $d=3$ in \cite{CLP12}, with the main results summarised in \cite[Chapter~5]{AC16}.
These results differ from the resolvent estimates discussed above, since they involve Sobolev norms of arbitrary order on spherical surfaces in $D_R$ (hence the radial derivative term in the $H^1(D_R)$ norm is not directly controlled, nor the $H^1(\Oin)$ norm). 
Some of these results describe in detail the behaviour of the solution when $\nin>1$, including the super-algebraic growth of the solution operator through a sequence of real wavenumbers, and we recap them in \S\ref{sec:blowup} below.
\ere

\bre\mythmname{Truncated transmission problems}
\label{rem:trunc}
When solving scattering problems on unbounded domains numerically, it is common to truncate the domain and impose a boundary condition to approximate the Sommerfeld radiation condition; the simplest such boundary condition is the impedance condition $\dn u - \ri k u=0$ (see, e.g., the discussion in \cite[\S5.1]{BaSpWu:16} and the references therein).
The truncated transmission problem is therefore equivalent to the interior impedance problem with piecewise-constant wavenumber. 
The paper \cite{CaVo:10} contains the analogue of the results in \cite{CaPoVo:99} for the truncated transmission problem (and in particular the bound \eqref{eq:res2} above).
Wavenumber-explicit bounds on this BVP have recently been obtained in \cite{Ch:16, BaChGo:16, SaTo:17, GrSa:18, GrPeSp:18} for real $n_i$ and \cite{OhVe:16} for complex $n_i$.
Apart from \cite{SaTo:17}, these recent investigations all use Morawetz identities (either explicitly or implicitly), with the impedance boundary condition dealt with as described in Remark \ref{rem:impbc} below (again, either explicitly or implicitly). 
The investigation \cite{SaTo:17} concerns the interior impedance problem in 1-d with piecewise-constant wavenumber, and uses 
the fact that,  in this case, the Green's function can be expressed in terms of the solution of a linear system.
\ere

\bre\mythmname{Transmission problems not involving a bounded obstacle}\label{rem:rough}
Identities related to those of Morawetz have also been used to prove results about (i) scattering by rough surfaces when the wave-number is piecewise constant and (ii) the transmission problem through an infinite penetrable layer (where in both cases the wavenumbers satisfy appropriate analogues of \eqref{eq:cond}). 
For (i) see \cite{ChZh:98, ZhCh:98, ZhCh:98a,HuLiQuZh:15}, and \cite[Chapter 2]{Th:06} 
(these works consider more  general classes of wave-number that include piecewise-constant cases), and for (ii) see \cite{ChZh:99, Fo:05,LeRi:10}, and \cite[Chapter 4]{Th:06}.
The identities used are essentially \eqref{eq:morid1} below with $\beta=0$ and the vector field $\bx$ replaced by a vector field perpendicular to the surface/layer.
\ere

\bre\mythmname{Transmission problems when $\Im \nin>0$}
\label{rem:complex}
Remark \ref{rem:complex2} below shows how analogues of Theorems  \ref{thm:FirstBound}--\ref{thm:nonzero} hold when $n_i \in \Com$ with $0<\Im n_i\leq \delta/k$ and $\delta$ is a sufficiently small constant (the occurrences of $n_i$ in the conditions \eqref{eq:cond}, \eqref{eq:condTracea}, and \eqref{eq:gDgNcond} are then replaced by $\Re n_i$).
This condition on the imaginary part is similar to that in \cite{NV12}, with this paper considering 
the Helmholtz transmission problem when $\Oi$ is the union of two concentric balls, modelling an inhomogeneity (with $n_i$ real) surrounded by an absorbing layer (with $n_i$ complex and $\Im n_i$ proportional to $1/k$).
Like the works \cite{AC16, Cap12, CLP12} discussed above, \cite{NV12} is interested in bounding the solution away from $\Oi$, but instead of using separation of variables, \cite{NV12} uses Morawetz identities to prove its bounds. 
The paper \cite{HPV07} proves bounds analogous to those in \cite{AC16, Cap12, CLP12} in the case when $\Oi$ is a ball and $\Im n_i>0$, again using separation of variables and bounds on Bessel and Hankel functions.
\ere

\section{Morawetz identities}\label{sec:Morawetz}

In this section we prove the identities that are the basis of the proofs of Theorems \ref{thm:FirstBound} and \ref{thm:nonzero}. The history of these identities is briefly discussed in Remark \ref{rem:biblio} below.

\ble\mythmname{Morawetz-type identity}\label{le:morid1}
Let $D\subset \Rea^d$, $d\geq 2$. Let $v\in C^2(D)$ and let $a,n,\alpha,$ $\beta\in\Rea$.
Let 
\beqs
\cL_{a,n} v:= a\Delta v + k^2 n\, v,
\eeqs
and let
\beq\label{eq:cM}
\cM v:= \bx\cdot \gv - \ri k \beta v + \alpha v.
\eeq
Then 
\begin{align}\nonumber 
2 \Re \big\{\overline{\cM v } \,\cL_{a,n} v \big\} = &\, \nabla \cdot \bigg[ 2 \Re\big\{\overline{\cM v}\, a \gv\big\} + \bx\big(k^2n \nvs -a\ngvs\big)\bigg] \\
&\hspace{10ex}
- (2\alpha -d +2) a\ngvs -(d-2\alpha)nk^2 \nvs.
\label{eq:morid1}
\end{align}
\ele

\bpf
This follows from expanding the divergence on the right-hand side of \eqref{eq:morid1}. Note that the identity \eqref{eq:morid1} is a special case of both \cite[Lemma 2.1]{SpKaSm:15} (where the multiplier $\cM u$ is generalised) and
\cite{GrPeSp:18} (where the operator $\cL_{a,n}$ is generalised).
\epf

The proofs of the main results are based on integrating the identity \eqref{eq:morid1} over $\Oin$ and $D_R(:= \Oout\cap B_R)$ and using the divergence theorem. Our next result, therefore, is an integrated version of \eqref{eq:morid1}. To state this result it is convenient to define the space
\beq\label{eq:V}
V(D):= \bigg\{ 
v\in H^1(D) :\;  \Delta v\in L^2(D), \;\dnu v \in L^2(\partial D),\; v \in H^1(\partial D)
\bigg\},
\eeq
where $D$ is a bounded Lipschitz open set with  outward-pointing unit normal vector $\bnu$.

\begin{lemma}
\mythmname{Integrated form of the Morawetz identity \eqref{eq:morid1}}
\label{lem:morid1int}
Let $D$ be a bounded Lipschitz open set, with boundary $\partial D$ and outward-pointing unit normal vector $\bnu$.
If $v \in V(D)$, $a,n,\alpha,\beta\in\Rea$, then
\begin{align}\nonumber 
&\int_D 2 \Re \big\{\overline{\cM v } \,\cL_{a,n} v \big\}
+ (2\alpha -d +2) a\ngvs  +(d-2\alpha) n k^2 \nvs
\\
&=\int_{\partial D}(\bx\cdot\bnu)\left(a \left|\dnu v\right|^2 -a|\nabla_T v|^2 + k^2 n \nvs\right)
+ 2\Re\Big\{\big(  \bx\cdot\overline{\nabla_T v}+ \ri k \beta \vb + \alpha \vb\big)a \dnu v\Big\}.
\label{eq:morid1int}
\end{align}
\ele
\begin{proof}
If $v\in C^\infty(\conj D)$, then \eqref{eq:morid1int} follows from divergence theorem
$\int_D \nabla \cdot \bF = \int_{\partial D} \bF \cdot \bnu$.
By \cite[Lemmas 2 and 3]{CoD98}, $C^\infty(\conj D)$ is dense in $V(D)$ and the result then follows 
since \eqref{eq:morid1int} is continuous in $v$ with respect to the topology of $V(D)$. 
\epf

The proofs of the main results use different multipliers in different domains. More precisely, we use
\begin{align*}
A_D A_N \bigg(& \bx\cdot \nabla u -\ri kR \sqrt{\frac{n_o}{a_o}} u + \frac{d-1}{2}u\bigg) \quad\tin \Oin,\\
&\bx\cdot \nabla u -\ri kR \sqrt{\frac{n_o}{a_o}} u + \frac{d-1}{2}u\,\,\,\,\quad\tin D_R,\tand\\
&\bx\cdot \nabla u -\ri kr \sqrt{\frac{n_o}{a_o}} u + \frac{d-1}{2}u\qquad\tin \Rea^d\setminus D_R,
\end{align*}
where $R$ is the radius of the ball in which we bound the solution, and $r:=|\bx|$.
The first two of these three multipliers are multiples of $\cM u$, as defined in \eqref{eq:cM}, with $\alpha=\frac{d-1}2$ and $\beta=R\sqrt{\nout/\aout}$; the third one is slightly different in that the coefficient $\beta=r\sqrt{\nout/\aout}$ depends on the position vector. Therefore, the identity arising from this last multiplier is not covered by Lemma \ref{le:morid1} but is given in the following lemma.

\begin{lemma}
\mythmname{Morawetz--Ludwig identity, \cite[Equation 1.2]{MoLu:68}}
\label{le:ML}
 Let $v \in C^2(D)$ for some $D\subset \Rea^d$, $d\geq 2$. 
Let $\kappa\in\IR$, $\cL v:=(\Delta +\kappa^2)v$ and let
 \beq\label{Malpha}
\cM_\alpha v :=	 r\left(v_r -\ri \kappa v + \frac{\alpha}{r}v\right),
\eeq
where $\alpha \in \Rea$ and $v_r=\bx\cdot \gv/r$. Then 
\bal\nonumber
2\Re\{ \overline{\cM_\alpha v} \cL v\} =  
&\,\nabla \cdot \bigg[2\Re \left\{\overline{\cM_{\alpha} v} \gv\right\}
+ \left(\kappa^2\nvs - \ngvs \right)\bx\bigg] \\&+ \big(2\alpha -(d-1)\big)\big(\kappa^2 \nvs - \ngvs\big) - \big(\ngvs -\nvrs\big)- \big| v_r -\ri \kappa v\big|^2.
\label{eq:ml2d}
\end{align}
\end{lemma}

The Morawetz--Ludwig identity \eqref{eq:ml2d} is a variant of the identity \eqref{eq:morid1} with $a=1$, $n=1$, $k=\kappa$, and $\beta=r$ (instead of being a constant); for a proof, see \cite{MoLu:68}, \cite[Proof of Lemma~2.2]{SpChGrSm:11}, or \cite[Proof of Lemma 2.3]{SpKaSm:15}.

As stated above, we use the Morawetz--Ludwig identity in $\Rea^d\setminus B_R$ (it turns out that this identity ``takes care'' of the contribution from infinity). It is convenient to encode the application of this identity in $\Rea^d\setminus B_R$ in the following lemma (slightly more general versions of which appear in \cite[Lemma 2.1]{CWM08} and 
\cite{GrPeSp:18}).

\ble\mythmname{Inequality on $\Gamma_R$ used to deal with the contribution from infinity} \label{lem:2.1}
Let $u$ be a solution of the homogeneous Helmholtz equation $\cL u=0$ in $\Rea^d\setminus \overline{B_{R_0}}$ (with $d\geq 2$), for some $R_0>0$, satisfying the Sommerfeld radiation condition. 
Then, for $R>R_0$, 
\beq\label{eq:2.1}
\int_{\Gamma_{R}} R\left( \left|\pdiff{u}{r}\right|^2 - |\nabla_{T} u|^2 + \kappa^2 |u|^2\right)   - 2 \kappa R\, \Im \int_{\Gamma_R} \bar{u} \pdiff{u}{r} + (d-1)\Re \int_{\Gamma_R}\bar{u}\pdiff{u}{r} \leq 0,
\eeq
where $\nabla_{T}$ is the tangential gradient on $r=R$ (recall that this is such that $\gv = \nabla_T v + \widehat{\bx}v_r$ on $r=R$).
\ele

\bpf
We now integrate \eqref{eq:ml2d} with $v=u$ and $2\alpha =d-1$ over $B_{R_1}\setminus B_R$, use the divergence theorem, and then let $R_1\tendi$
(note that using the divergence theorem is allowed since $u$ is $C^\infty$ by elliptic regularity).

Writing the  identity  \eqref{eq:ml2d} as $\nabla\cdot\bQ(v)=P(v)$, we have that 
if $u$ is a solution of $\cL u=0$ in $\Rea^d\setminus \overline{B_{R_1}}$ satisfying the Sommerfeld radiation condition \eqref{eq:src}, then
\beqs
\int_{\Gamma_{R_1}} \bQ(u) \cdot\widehat{\bx}
 \tendo \quad \tas R_1\tendi
\eeqs
(independent of the value of $\alpha$ in the multiplier $\cM_\alpha u$); see \cite[Proof of Lemma 5]{MoLu:68}, \cite[Lemma~2.4]{SpChGrSm:11}. 
Note that, although the set-up of \cite{SpChGrSm:11} is for $d=2, 3$, the proof of \cite[Lemma~2.4]{SpChGrSm:11} holds for $d\geq 2$.

Then, using the decomposition $\gv = \nabla_T v + \widehat{\bx}v_r$ on the integral over $\Gamma_R$ (or equivalently the right-hand side of \eqref{eq:morid1int}),
we obtain that
\begin{align*}
\int_{\Gamma_R} \bQ(u) \cdot\widehat{\bx} &=
\int_{\Gamma_R} R\left( \left|\pdiff{u}{r}\right|^2  - |\nabla_T u|^2+ \kappa^2 |u|^2\right)   - 2 \kappa R\, \Im \int_{\Gamma_R} \bar{u} \pdiff{u}{r} + (d-1)\Re \int_{\Gamma_R}\bar{u}\pdiff{u}{r}  \nonumber \\
& =  -\int_{\Rea^d\setminus B_R}\left(\big(\ngus -\nurs\big)+ \left|
 u_r - \ri \kappa u \right|^2\right)\leq 0.
\end{align*}
\epf

\bre\mythmname{Far-field impedance boundary condition}\label{rem:impbc}
If the infinite domain $\Oout$ is truncated, and the radiation condition approximated by an impedance boundary condition, then 
an analogous inequality to that in Lemma \ref{lem:2.1} holds; see \cite[Lemma~4.6]{GrPeSp:18}. 
This analogous inequality allows one to extend the results of Theorems \ref{thm:FirstBound} and \ref{thm:nonzero} to this truncated BVP (as mentioned in Remark \ref{rem:extensions} above); see \cite{GrPeSp:18}  for more details.
\ere

\bre\mythmname{Bibliographic remarks}\label{rem:biblio}
The multiplier $\bx \cdot \gv$ was introduced by Rellich in \cite{Re:40}, and has been well-used since then in the study of the Laplace, Helmholtz, and other elliptic equations, see, e.g., the references in \cite[\S5.3]{CGLS12}, \cite[\S1.4]{MOS12}. 

The idea of using a multiplier that is a linear combination of derivatives of $v$ and $v$ itself, such as $\cM v$, 
is attributed by Morawetz in \cite{Mo:61} to Friedrichs. The multiplier $\cM_\alpha v$ \eqref{Malpha} for the Helmholtz equation was introduced by Morawetz and Ludwig in \cite{MoLu:68} and the multiplier $\cM v$ \eqref{eq:cM} (with $\bx$ replaced by a general vector field and $\alpha$ and $\beta$ replaced by general scalar fields) is implicit in Morawetz's paper \cite{Mo:75} (for more discussion on this, see \cite[Remark 2.7]{SpKaSm:15}).
\ere

\section{Proofs of the main results}\label{sec:proofs}

\bpf[Proof of Theorem \ref{thm:FirstBound}]
We use the integrated Morawetz identity \eqref{eq:morid1int} with first $D=\Oi$ and then $D=D_R$. In both cases we take $2\alpha=d-1$ and use the same (as yet unspecified) $\beta$; in the first case we take
$v=u_i$, $a=a_i, n=n_i$, and in the second case we take $v=u_o$, $a=a_o, n=n_o$.
Using \eqref{eq:morid1int} is justified since the regularity result in Lemma \ref{lem:exist} shows that $u_i \in V(\Oi)$ and $u_o\in V(D_R)$.
We get 
\begin{align}\nonumber 
&\int_{\Oi} a_i|\nabla u_i|^2 + n_i k^2 |u_i|^2 \\
&\qquad=- 2 \Re \int_{\Oi} \overline{\cM u_i }\, f_i
+\int_{\Gamma}(\bx\cdot\bn)\left(a_i \left|\dn u_i\right|^2 -a_i|\nabla_T u_i|^2 + k^2 n_i |u_i|^2\right)
\nonumber\\&\hspace{43mm}
+ 2\Re\left\{\left(  \bx\cdot\overline{\nabla_T u_i}+ \ri k \beta\overline{u_i}  + \frac{d-1}{2} \overline{u_i}\right)a_i \dn u_i\right\}\label{eq:1}
\end{align}
and 
\begin{align}\nonumber 
&\int_{D_R} a_o|\nabla u_o|^2 + n_o k^2 |u_o|^2 \\
&=- 2 \Re \int_{D_R} \overline{\cM u_o}\, f_o
-\int_{\Gamma}(\bx\cdot\bn)\left(a_o \left|\dn u_o\right|^2 -a_o|\nabla_T u_o|^2 + k^2 n_o|u_o|^2\right)
\nonumber\\&\hspace{43mm}
+ 2\Re\left\{\left(  \bx\cdot\overline{\nabla_T u_o}+ \ri k \beta\overline{u_o}  + \frac{d-1}{2} \overline{u_o}\right)a_o\dn u_o\right\}
\nonumber\\
&\hspace{3mm}+\int_{\Gamma_R} R \left( a_o \left|\pdiff{u_o}{r}\right|^2 - a_o |\nabla_T u_o|^2 + k^2 n_o |u_o|^2 \right) 
\nonumber\\&\hspace{3mm}
- 2a_o k \beta \Im \int_{\Gamma_R} \overline{u_o}\pdiff{u_o}{r} + a_o(d-1) \Re \int_{\Gamma_R} \overline{u_o}\pdiff{u_o}{r}.\label{eq:2}
\end{align}
Multiplying the inequality \eqref{eq:2.1} by $a_o$ and letting $\kappa= k \sqrt{n_o/a_o}$, we see that if we choose $\beta= R\sqrt{n_o/a_o}$ then the terms on $\Gamma_R$ on the right-hand side of \eqref{eq:2} are non-positive.

We then multiply \eqref{eq:2} by an arbitrary $\eta>0$ and add to \eqref{eq:1} to get
\begin{align}\nonumber
&\int_{\Omega_i}(a_i|\nabla u_i|^2+k^2n_i|u_i|^2)
+\eta\int_{D_R}(a_o|\nabla u_o|^2+k^2n_o|u_o|^2)\\ \nonumber
&\le- 2 \Re \int_{\Oi} \overline{\cM u_i }\, f_i- 2\eta \Re \int_{D_R} \overline{\cM u_o}\, f_o\\
& + \int_\Gamma (\bx\cdot\bn)\left(a_i \left|\dn u_i\right|^2-\eta a_o \left|\dn u_o\right|^2 
-a_i|\nabla_T u_i|^2 +\eta a_o|\nabla_T u_o|^2
 + k^2 n_i|u_i|^2-\eta k^2 n_o|u_o|^2\right)\nonumber\\
&\hspace{10mm}+ 2\Re\left\{\left(  \bx\cdot\overline{\nabla_T u_i}+ \ri k \sqrt{\frac{n_o}{a_o}}R\overline{u_i}  + \frac{d-1}{2} \overline{u_i}\right)a_i\dn u_i\right\}\nonumber\\
&\hspace{10mm}-2 \eta \Re\left\{\left(  \bx\cdot\overline{\nabla_T u_o}+ \ri k \sqrt{\frac{n_o}{a_o}}R\overline{u_o}  + \frac{d-1}{2} \overline{u_o}\right)a_o\dn u_o\right\}.
\label{eq:3}
\end{align}
The volume terms at the right-hand side are bounded above by 
\begin{align}\nonumber
&\N{f_i}_{\Omega_i}\left( 2 \diam(\Omega_i) \N{\nabla u_i}_{\Omega_i} + \left(2 k \sqrt{\frac{n_o}{a_o}}R + d-1\right)\N{u_i}_{\Omega_i}\right)\\
&\qquad+\eta\N{f_o}_{D_R}\left(2R\N{\nabla u_o}_{D_R} +  \left(2 k \sqrt{\frac{n_o}{a_o}}R + d-1\right)\N{u_o}_{\Omega_o}\right)\label{eq:4}
\end{align}
(recall that $\|\cdot\|_{\Omega_{i/o}}$ denotes the $L^2$ norm on $\Omega_{i/o}$).
We now focus on the terms on $\Gamma$, and recall that we are assuming that $g_D=g_N=0$. 
Our goal is to choose $\eta$ so that the terms without a sign (i.e.~those on the last two lines of \eqref{eq:3}) cancel. Using the transmission conditions in \eqref{eq:BVP}, we see that this cancellation occurs if $\eta = 1/(A_D A_N)$. (It is at this point that we need $A_D$ and $A_N$ to be real; indeed, if the product $A_D A_N$ has non-zero real and imaginary parts, we cannot chose even a complex $\eta$ to cancel these terms).

Making this choice of $\eta$ and using the transmission conditions, we see that
the remaining terms on $\Gamma$ become
\beq
\int_\Gamma (\bx\cdot\bn) \left(
k^2 n_i |u_i|^2 \left (1 - \frac{A_D}{A_N}\frac{n_o}{n_i}\right)+
a_i |\dn u_i|^2 \left (1 - \frac{A_N}{A_D}\frac{a_i}{a_o}\right) 
+a_i |\nabla_T u_i |^2 
\left (-1 +\frac{A_D}{A_N}\frac{a_o}{a_i}\right) 
\right).
\label{eq:GammaTerm}
\eeq
These terms are negative and thus can be neglected if
\beqs
A_N n_i \leq A_D n_o\quad\tand \quad A_N a_i \geq A_D a_o,
\eeqs
or equivalently if \eqref{eq:cond} holds.

In summary, under the conditions \eqref{eq:cond}, we have that 
\beqs
a_i\N{\nabla u_i}_{\Oin}^2+k^2n_i\N{u_i}^2_{\Oin}
+\frac1{A_D A_N}\left(a_o\N{\nabla u_o}_{D_R}^2+k^2n_o\N{u_o}^2_{D_R}\right)
\eeqs
is bounded by \eqref{eq:4}. 
Using the Cauchy--Schwarz and Young inequalities we obtain the assertion of Theorem~\ref{thm:FirstBound}.
\epf

\bre\label{rem:unnatural}\mythmname{The origin of the condition \eqref{eq:cond}}
Condition \eqref{eq:cond} comes from requiring that each of the terms in \eqref{eq:GammaTerm} are non-positive. These terms are not independent, however, since they all depend on $u_i$. Despite this connection, we have not been able to lessen the requirements of \eqref{eq:cond} using only these elementary arguments, other than in Proposition \ref{thm:BoundViaTrace2} below where the term on $\Gamma$ involving $|u_i|^2$ is controlled by the norms of $u_i$ in $\Oin$ via a trace inequality.
\ere

\bpf[Proof of Theorem \ref{thm:nonzero}]
In the proof of Theorem~\ref{thm:FirstBound}, the assumption $g_D=g_N=0$ was used to derive
\eqref{eq:GammaTerm} from \eqref{eq:3}.
Now that $g_D$ and $g_N$ are not necessarily zero, we expand the terms on $\Gamma$ appearing in \eqref{eq:3} using the transmission conditions with $g_D,g_N\ne0$ and the fact that $\eta=1/\AD\AN$.
We control these terms on $\Gamma$ using $|\bx|\le\diam(\Oin)$ and $-\bx\cdot\bn\le-\gamma\diam(\Oin)$ for a.e.\ $\bx\in\Gamma$, and $(\ain\AN-\aout\AD)>0$ and $(\nout\AD-\nin\AN)>0$ from assumption~\eqref{eq:gDgNcond}.
We apply the weighted Young's inequality to nine terms, denoted $T_1,\ldots,T_9$, using positive coefficients $\xi_1,\ldots,\xi_9$:
\begin{align*}
&\int_\Gamma(\bx\cdot\bn)\bigg(a_i \left|\dn u_i\right|^2 -a_i|\nabla_T u_i|^2 
 + k^2 n_i|u_i|^2
 \\
&-\frac1{\aout\AD\AN} \left|\AN \ain\dn u_i+g_N\right|^2  
+ \frac\aout{\AD\AN}|\AD\nabla_T u_i+\nabla_T\gD|^2
 - k^2 \frac\nout{\AD\AN}|\AD u_i+\gD|^2\bigg)\\
&+ 2\Re\left\{\left(  \bx\cdot\overline{\nabla_T u_i}+ \ri k \sqrt{\frac{n_o}{a_o}}R\overline{u_i}  + \frac{d-1}{2} \overline{u_i}\right)a_i\dn u_i\right\}
\\
&-2  \Re\left\{\left(  \bx\cdot\Big(\overline{\nabla_T u_i+\frac1\AD\nabla_T\gD}\Big)
+\Big( \ri k \sqrt{\frac{n_o}{a_o}}R+\frac{d-1}{2}\Big)\overline{\Big(u_i+\frac1\AD\gD\Big)} \right)
\Big(a_i\dn u_i+\frac1\AN g_N\Big)\right\}
\hspace{-5mm}
\\
%%%%%%%%%%%%%%%%%%%%%%%%%%%%%%%%%%%%%%%
=&\int_\Gamma
(\bx\cdot\bn)\bigg(\ain\Big(1-\frac{\ain\AN}{\aout\AD}\Big) |\dn u_i|^2
+\ain\Big(\frac{\aout\AD}{\ain\AN}-1\Big)|\nabla_T u_i|^2
+k^2\nin\Big(1-\frac{n_o\AD}{\nin\AN}\Big)|u_i|^2
\\&
-\frac1{\aout\AD\AN} |g_N|^2
+\frac\aout{\AD\AN}|\nabla_T\gD|^2
-k^2 \frac\nout{\AD\AN}|\gD|^2
\\&
-2 \frac{a_i}{\aout\AD}\Re\underbrace{\{\dn u_i \conj{g_N}\}}_{T_1}
+2\frac\aout\AN\Re\underbrace{\{\nabla_T u_i\cdot\nabla_T\conj\gD\}}_{T_2}
-2k^2\frac\nout{\AN}\Re\underbrace{\{u_i\conj\gD\}}_{T_3}\bigg)
\\&
-2\frac1{\AN}\Re\underbrace{\{\bx\cdot\nabla_T \conj{u_i}g_N\}}_{T_4}
-2\frac1{\AN}\Re\underbrace{\Big\{\Big( \ri k \sqrt{\frac{n_o}{a_o}}R+\frac{d-1}{2}\Big)\conj u_ig_N\Big\}}_{T_5}
\\&
-2\frac{\ain}{A_D}\Re\underbrace{\left\{\bx\cdot\nabla_T\conj{\gD} \dn u_i\right\}}_{T_6}
-2\frac{\ain}{A_D}\Re\underbrace{\left\{\Big( \ri k \sqrt{\frac{n_o}{a_o}}R+\frac{d-1}{2}\Big)\overline{\gD}
\dn u_i\right\}}_{T_7}
\\&
-2\frac1{\AN\AD}\Re\underbrace{\left\{\bx\cdot\nabla_T\conj{\gD}  g_N\right\}}_{T_8}
-2\frac1{\AN\AD}\Re\underbrace{\left\{\Big( \ri k \sqrt{\frac{n_o}{a_o}}R+\frac{d-1}{2}\Big)\overline{\gD}
 g_N\right\}}_{T_9}
\\
%%%%%%%%%%%%%%%%%%%%%%%%%%%%%%%%%%%%%%%
\le&
-\gamma\diam(\Oin)\bigg(\ain\frac{\ain\AN-\aout\AD}{\aout\AD}\N{\dn u_i}_\Gamma^2
\\&\hspace{20mm}
+\frac{\ain\AN-\aout\AD}{\AN}\N{\nabla_T u_i}_\Gamma^2
+k^2\frac{\nout\AD-\nin\AN}{\AN}\N{u_i}_\Gamma^2\!\bigg)
\\
&-\gamma\diam(\Oin)\frac1{\aout\AD\AN} \N{g_N}_\Gamma^2
+\diam(\Oin)\frac\aout{\AD\AN}\N{\nabla_T\gD}_\Gamma^2
-k^2 \gamma\diam(\Oin)\frac\nout{\AD\AN}\N{\gD}_\Gamma^2
\\
&+\xi_1\ain\gamma\diam(\Oin)\frac{\ain\AN-\aout\AD}{\aout\AD}\N{\dn u_i}_\Gamma^2
+\frac{\ain\diam(\Oin)}{\xi_1\gamma\aout\AD(\ain\AN-\aout\AD)}\N{g_N}_\Gamma^2
\\
&+\xi_2\gamma\diam(\Oin)\frac{\ain\AN-\aout\AD}{\AN}\N{\nabla_T u_i}_\Gamma^2
+\frac{\diam(\Oin)\aout^2}{\xi_2\gamma\AN(\ain\AN-\aout\AD)}\N{\nabla_T\gD}_\Gamma^2
\\
&+\xi_3k^2\gamma\diam(\Oin)\frac{\nout\AD-\nin\AN}{\AN}\N{u_i}_\Gamma^2
+\frac{k^2\diam(\Oin)n_o^2}{\xi_3\gamma\AN(\nout\AD-\nin\AN)}\N{\gD}_\Gamma^2
\\
&+\xi_4\gamma\diam(\Oin)\frac{\ain\AN-\aout\AD}{\AN}\N{\nabla_T u_i}_\Gamma^2
+\frac{\diam(\Oin)}{\xi_4\gamma\AN(\ain\AN-\aout\AD)}\N{g_N}_\Gamma^2
\\
&+\xi_5k^2\gamma\diam(\Oin)\frac{\nout\AD-\nin\AN}{\AN}\N{u_i}_\Gamma^2
+\frac{\Big|\ri k \sqrt{\frac{n_o}{a_o}}R+\frac{d-1}{2}\Big|^2}{\xi_5\gamma\diam(\Oin) k^2\AN(\nout\AD-\nin\AN)}
\N{g_N}_\Gamma^2
\\
&+\xi_6\ain\gamma\diam(\Oin)\frac{\ain\AN-\aout\AD}{\aout\AD}\N{\dn u_i}_\Gamma^2
+\frac{\ain\aout\diam(\Oin)}{\xi_6\gamma\AD(\ain\AN-\aout\AD)}\N{\nabla_T\gD}_\Gamma^2
\\
&+\xi_7\ain\gamma\diam(\Oin)\frac{\ain\AN-\aout\AD}{\aout\AD}\N{\dn u_i}_\Gamma^2
+\frac{\ain\aout\Big|\ri k \sqrt{\frac{n_o}{a_o}}R+\frac{d-1}{2}\Big|^2}{\xi_7\gamma\AD\diam(\Oin)(\ain\AN-\aout\AD)}\N{\gD}_\Gamma^2
\\
&+\xi_8\frac{\diam(\Oin)}{\aout\AN\AD}\N{g_N}_\Gamma^2
+\frac{\diam(\Oin)\aout}{\xi_8\AN\AD}\N{\nabla_T\gD}_\Gamma^2
\\&
+\xi_9\frac{\diam(\Oin)}{\aout\AN\AD} \N{g_N}_\Gamma^2
+\frac{\aout\Big|\ri k \sqrt{\frac{n_o}{a_o}}R+\frac{d-1}{2}\Big|^2}
{\xi_9\diam(\Oin)\AD\AN} \N{\gD}_\Gamma^2.
\end{align*}
We choose the weights as
$$\xi_1=\xi_6=\xi_7=\frac13,\quad \xi_2=\xi_3=\xi_4=\xi_5=\frac12,\quad \xi_8=\xi_9=1,$$
so that all terms containing $u_i$ cancel each other, and we are left with
\begin{align*}
%%%%%%%%%%%%%%%%%%%%%%%%%%%%%%%%%%%%%%%%%%%%%%%%%%%%%%%%%
&\frac{\diam(\Oin)}{\AN}\bigg[
-\frac{\gamma\nout}\AD+\frac{2n_o^2}{\gamma(\nout\AD-\nin\AN)}
\\&\hspace{20mm}
+\frac{1}{\AD\diam(\Oin)^2}\Big(1+\frac{3\ain\AN}{\gamma(\ain\AN-\aout\AD)}\Big)
\Big(\nout R^2+\frac{\aout(d-1)^2}{4k^2}\Big)
\bigg]k^2\N{\gD}_\Gamma^2
\\
&+\frac{\diam(\Oin)\aout}{\AD\AN}\bigg[
2+\frac{2\aout\AD+3\ain\AN}{\gamma(\ain\AN-\aout\AD)}\bigg]\N{\nabla_T\gD}_\Gamma^2
\\
&+\frac{\diam(\Oin)}{\aout\AN\AD}\bigg[
2-\gamma+\frac{3\ain\AN+2\aout\AD}{\gamma(\ain\AN-\aout\AD)}
+\frac{2\AD\Big(\nout R^2+\frac{\aout(d-1)^2}{4k^2}\Big)}{\gamma\diam(\Oin)^2(\nout\AD-\nin\AN)}
\bigg]\N{g_N}_\Gamma^2
\\
\le&
%%%%%%%%%%%%%%%%%%%%%%%%%%%%%%%%%%%%%%%%%%%%%%%%%%%%%%%%%
\bigg[
\frac{2\diam(\Oin)n_o^2}{\gamma\AN(\nout\AD-\nin\AN)}
+\frac{(3+\gamma)\ain\Big(\nout R^2+\frac{\aout(d-1)^2}{4k^2}\Big)}{\AD\gamma\diam(\Oin)(\ain\AN-\aout\AD)}
\bigg]k^2\N{\gD}_\Gamma^2
\\
&+\frac{\diam(\Oin)\aout\big((3+2\gamma) \ain \AN + 2 \aout\AD\big)}{\AD\AN\gamma(\ain\AN-\aout\AD)}\N{\nabla_T\gD}_\Gamma^2
\\
&+\frac{1}{\gamma\aout\AN\AD}\bigg[
\frac{\diam(\Oin)(4\ain\AN+2\aout\AD)}{\ain\AN-\aout\AD}
+\frac{2\AD\Big(\nout R^2+\frac{\aout(d-1)^2}{4k^2}\Big)}{\diam(\Oin)(\nout\AD-\nin\AN)}
\bigg]\N{g_N}_\Gamma^2,
\end{align*}
where we used also $0<\gamma\le1/2$ and dropped some negative terms.
Then the bound in the assertion, \eqref{eq:gDgNbound}, follows as in the proof of Theorem \ref{thm:FirstBound}, recalling that the use of Young's inequality for the volume norms gives a further factor of 2 in front of norms on $\Gamma$ in the final bound.
\epf

\bre\label{rem:complex2}\mythmname{Extensions of Theorems \ref{thm:FirstBound}--\ref{thm:nonzero} to the case when $\Im n_i>0$}
We now explain how analogues of Theorems~\ref{thm:FirstBound}--\ref{thm:nonzero} hold when $n_i\in \Com$ with $0<\Im n_i\leq \delta/k$ and $\delta$ is sufficiently small (the occurrences of $n_i$ in the conditions \eqref{eq:cond}, \eqref{eq:condTracea}, and \eqref{eq:gDgNcond} are then replaced by $\Re n_i$).
We first consider Theorem \ref{thm:FirstBound}. 
Under the assumption that the existence, uniqueness, and regularity results of Lemma \ref{lem:exist} hold for the boundary value problem with such $n_i$,
Equation \eqref{eq:3} now holds with $n_i$ replaced by $\Re n_i$ and $f_i$ replaced by $f_i -\ri k^2 (\Im n_i) u_i$
We therefore have the extra term 
\beqs
-2\Im\int_{\Oin} \overline{\cM u_i} \,k^2 (\Im n_i) u_i 
\eeqs
on the right-hand side of \eqref{eq:3}. If $0<\Im n_i\leq \delta/k$ and $\delta$ is sufficiently small, then this term can be absorbed into the weighted $H^1$-norm of $u_i$ on the left-hand side of \eqref{eq:3} (using the Cauchy--Schwarz and weighted Young inequalities). The result is a bound with the same $k$-dependence as \eqref{eq:bound1}, 
but slightly different constants on the right-hand side. This proof of an a priori bound, under the assumption of existence, implies uniqueness. One can then check that the proof of existence and regularity in Appendix \ref{sec:appA} goes through with $n_i$ of this particular form. 
Finally, the extensions to the proof of Theorem \ref{thm:FirstBound} needed to prove Theorem \ref{thm:nonzero} go through as before (since these only involve the terms on $\Gamma$).
\ere

\bpf[Proof of Theorem \ref{thm:res}]

The result \cite[Lemma 2.3]{Vo:99} implies that the assertion of the theorem will hold if 
(i)  $\Rimp(k)$ is holomorphic for  $\Im k> 0$ and (ii) there exist $C_4>0$ and $k_0>0$ such that
\beq\label{eq:pert1}
\N{\Rimp(k)}_{{L^2(\Oin)\oplus L^2(\Oout) \rightarrow L^2(\Oin)\oplus L^2(\Oout)}}\leq \frac{C_4}{k} \quad \tfa k\geq k_0.
\eeq
Note that, in applying Vodev's result, we take Vodev's obstacle $\Omega$ to be the empty set, $N_0=1$, $\Omega_1$ equal to our $\Omega_i$, and $g_{ij}^{(1)}=\delta_{ij}$. 
We also note that the set-up in \cite{Vo:99} assumes that $\Oi$ is smooth. Nevertheless, the result \cite[Lemma 2.3]{Vo:99} boils down to a perturbation argument (via Neumann series) and a result about the free resolvent (i.e.~the inverse of the Helmholtz operator in the absence of any obstacle) \cite[Lemma~2.2]{Vo:99}; both of these results  are independent of $\Oi$, and so \cite[Lem\-ma~2.3]{Vo:99} is valid when $\Oi$ is Lipschitz.

Since $\Oi$ is star-shaped and the condition \eqref{eq:cond} is satisfied, Theorem \ref{thm:FirstBound} implies that the estimate on the real axis \eqref{eq:pert1} holds, and thus we need only show that 
$\Rimp(k)$ is holomorphic on $\Im k> 0$.

Observe that $\Rimp(k)$ is well-defined for $\Im k\geq 0$ by the existence and uniqueness results of Lemma \ref{lem:exist}.
Analyticity follows by applying the Cauchy--Riemann
operator $\partial/\partial \overline{k}$ to the BVP \eqref{eq:BVP}. Indeed, 
by using Green's integral representation in $\Omega_{i/o}$ we find that $\partial (\Delta u)/\partial \overline{k} = \Delta( \partial u/\partial \overline{k})$, and similarly for $\nabla u$. These in turn imply that $\dn (\partial u/\partial \overline{k})= \partial (\dn u)/\partial \overline{k}$. Therefore, applying $\partial/\partial \overline{k}$ to \eqref{eq:BVP}, we find that $\partial u/\partial \overline{k}$ satisfies the Helmholtz transmission problem with zero volume and boundary data, and thus must vanish by the uniqueness result.
\epf

The condition \eqref{eq:cond} in Theorem \ref{thm:FirstBound} implies that $\nin/\ain\le\nout/\aout$, namely that the wave\-length $\lambda=(2\pi\sqrt a)/(\sqrt nk)$ of the solution $u$ is larger in the inner domain $\Oin$ than in $\Oout$.
In the next proposition we extend the result of Theorem \ref{thm:FirstBound} to a case where this condition is slightly violated.

\begin{proposition}\label{thm:BoundViaTrace2}
Assume that $\Oi$ is star-shaped,
\beq\label{eq:condTracea}
\frac{\nout}{\nin}
\left(\frac{\nin}{\nout}- \frac{\AD}{\AN}\right)
\left(d+ \sqrt{d^2 + \frac{4n_i}{a_i}\big(k\diam(\Oin)\big)^2}\right)
<1,
\qquad
\frac{A_D}{A_N}\leq \frac{a_i}{a_o},
\eeq
$g_N=g_D=0$  and $k>0$.
Then the solution of BVP \eqref{eq:BVP} satisfies
\begin{align}
&G\left(a_i\N{\nabla u_i}_{\Oin}^2+k^2n_i\N{u_i}^2_{\Oin}\right)
+\frac1{A_D A_N}\left(a_o\N{\nabla u_o}_{D_R}^2+k^2n_o\N{u_o}^2_{D_R}\right)
\nonumber\\
&\hspace{2cm}\le
\bigg[\frac{4\diam(\Omega_i)^2}{a_i}+\frac{1}{n_i}\left(2\sqrt{\frac{n_o}{a_o}} R +\frac{d-1}{k}\right)^2\bigg]
\N{f_i}_{\Omega_i}^2
\nonumber\\
&\hspace{3cm}+\frac1{A_D A_N}\bigg[\frac{4R^2}{a_o}+\frac{1}{n_o}\left(2\sqrt{\frac{n_o}{a_o}} R +\frac{d-1}{k}\right)^2\bigg]\N{f_o}_{D_R}^2,
\label{eq:bound2a}
\end{align}
where $G$ is defined by
\beq\label{eq:G2}
G:=\half \left(1-
\frac{\nout}{\nin}
\left(\frac{\nin}{\nout}- \frac{\AD}{\AN}\right)
\left(d+ \sqrt{d^2 + \frac{4n_i}{a_i}\big(k\diam(\Oin)\big)^2}\right)
\right)
\eeq
and is positive by the first inequality in \eqref{eq:condTracea}.
\end{proposition}

To better understand the condition \eqref{eq:condTracea}, consider the simple case when $\nin=\ain=\aout=\AD=\AN=1$. Then the condition \eqref{eq:condTracea} is satisfied if 
\beq\label{eq:condTrace2a}
\nout>1-\frac1{d + \sqrt{d^2+4(k\diam(\Oi))^2}}.
\eeq
For a fixed $\nout>1$ and sufficiently small, the condition \eqref{eq:condTrace2a} is an upper bound on $k$ under which the estimate \eqref{eq:bound2a} holds---this is consistent with the results on super-algebraic growth in $k$ of the resolvent for $\nout>1$ recapped in \S\ref{sec:blowup} below.
If we allow $\nout$ to be a function of $k$, then the condition \eqref{eq:condTrace2a} implies that the estimate \eqref{eq:bound2a} holds if the distance between $\nout$ and $\nin(=1)$ decreases like $1/k$ as $k\tendi$.

\bpf[Proof of Proposition \ref{thm:BoundViaTrace2}]
The proof proceeds exactly the same was as the proof of Theorem \ref{thm:FirstBound} up to \eqref{eq:GammaTerm}. Now, the assumption in \eqref{eq:condTracea} that $A_D a_o \leq A_N a_i$ implies that the terms in \eqref{eq:GammaTerm} are bounded by 
\beq\label{eq:GammaTerm2}
\int_\Gamma (\bx\cdot\bn) \left(
k^2 n_i |u_i|^2 \left (1 - \frac{A_D}{A_N}\frac{n_o}{n_i}\right)\right).
\eeq
For all $w\in H^1(\Oin)$ and $\epsilon>0$, we have the following weighted trace inequality:
\begin{align}\nonumber
\int_\Gamma (\bx\cdot\bn)|w|^2
&=\int_\Oin \div(\bx|w|^2)\\&
=\int_\Oin (d|w|^2+2\bx\cdot\Re\{w\nabla \conj w\})
\le (d+\epsilon)\N{w}^2_{\Oin}+\frac1\epsilon\diam(\Oin)^2\N{\nabla w}^2_{\Oin}.
\label{eq:trace}
\end{align}
We chose $\epsilon$ so that 
\beq\label{eq:euan2}
\frac{d+\epsilon}{n_i k^2} =\frac{1}{\epsilon} \frac{\diam(\Oi)^2}{a_i},
\eeq
so that the right-hand side of \eqref{eq:trace} becomes
\beqs
\frac{d+\epsilon}{n_i k^2}\left( n_i k^2 \N{w}^2_{\Oin} + a_i \N{\nabla w}^2_{\Oin}\right).
\eeqs
The requirement \eqref{eq:euan2}, and the fact that $\epsilon>0$, imply that
\beqs
\epsilon = \half \left( - d + \sqrt{d^2 + \frac{4n_i}{a_i} \big(k \diam(\Oi)\big)^2}\right).
\eeqs
We then get that \eqref{eq:GammaTerm2} is bounded by
\beqs
 \frac{n_o}{n_i}\left(\frac{\nin}{\nout}- \frac{\AD}{\AN}\right)
\half\left(d+ \sqrt{d^2 + \frac{4n_i}{a_i}\big(k\diam(\Oin)\big)^2}\right)\big( n_i k^2\N{u_i}^2_{\Oin} + a_i \N{\nabla u_i}^2_{\Oin}\big).
\eeqs
The requirement \eqref{eq:condTracea} implies that this term is strictly less than 
$( n_i k^2\N{u_i}^2_{\Oin} + a_i \N{\nabla u_i}^2_{\Oin})/2$,
and thus the argument proceeds as before (with the other half of the norm being used to deal with the terms in \eqref{eq:4} via the Cauchy--Schwarz and Young inequalities).
\epf

\section{Super-algebraic growth in \texorpdfstring{$k$}{k} when the condition \texorpdfstring{\eqref{eq:cond}}{} is violated}\label{sec:blowup}
In \S\ref{sec:PoVo} we adapt the results of \cite{PoVo:99}, for $C^\infty$ convex $\Oi$, on super-algebraic growth of the solution operator through a sequence of complex wavenumbers when the condition \eqref{eq:cond} does not hold to prove an analogous result for real wavenumbers. In \S\ref{sec:ball} we briefly highlight the analogous results of \cite{Cap12, CLP12, AC16}, valid for real wavenumbers, when $\Oi$ is a ball.
As explained at the end of \S\ref{sec:intro}, our main motivation for doing this is the recent interest in \cite{Ch:16, BaChGo:16, OhVe:16, SaTo:17, GrSa:18} on 
how the solution of the interior impedance problem with piecewise-constant wavenumber (which is an approximation of the transmission problem) depends on the wavenumber, and the fact that the results in this section partially answer questions/conjectures from \cite{BaChGo:16, SaTo:17}.

\subsection{Adapting the results of Popov and Vodev \texorpdfstring{\cite{PoVo:99}}{}}\label{sec:PoVo}

The paper \cite{PoVo:99} considers the Helmholtz transmission problem \eqref{eq:BVP} with 
$\Oin$ a $C^\infty$ convex domain with strictly positive curvature, $\ain=\aout=\nout=\AD=1$ and $g_D=g_N=0$; i.e., the BVP \eqref{eq:BVP2}.

From our point of view, the importance of \cite{PoVo:99}
is that their results can be adapted to show that if $\Oin$ is a $C^\infty$ convex domain with strictly positive curvature, $n_i>1$ and $A_N>0$, then there exists an increasing sequence of real wavenumbers through which the solution operator grows super-algebraically; this is stated as Corollary \ref{cor:PoVo2} below. To state the results of \cite{PoVo:99} we need to define a \emph{quasimode}.

\begin{defin}\mythmname{Quasimode for the BVP \eqref{eq:BVP2}}\label{def:quasi}
A quasimode is a sequence
\beqs
\big\{ k_j , \big(u_i^{(j)}, u_o^{(j)}\big) \big\}_{j=1}^\infty
\eeqs
where $k_j \in \Com$, $|k_j|\tendi$, $\Re k_j \geq 1$, $u_{i/o}^{(j)} \in C^\infty(\overline{\Omega_{i/o}})$, the support of $u^{(j)}_{i/o}$ is contained in a fixed compact neighbourhood of $\Gamma$, $\|u_i^{(j)}\|_{L^2(\Gamma)}=1$, 
\begin{subequations}
\begin{align}
\N{ (\Delta +n_i k_j^2)u_i^{(j)}}_{L^2(\Oin)} &= \cO(|k_j|^{-\infty}),\\
\N{ (\Delta +k_j^2)u_o^{(j)}}_{L^2(\Oout)} &= \cO(|k_j|^{-\infty}),\\
\N{ u_i^{(j)}-u_o^{(j)} }_{H^2(\Gamma)} &= \cO(|k_j|^{-\infty}),\quad \tand \label{eq:QM3}\\
\N{ \dn u_i^{(j)}-A_N \dn u_o^{(j)} }_{H^2(\Gamma)} &= \cO(|k_j|^{-\infty}),\label{eq:QM4}
\end{align}
\end{subequations}
where, given an infinite sequence of complex numbers $\{z_j\}_{j=1}^\infty$, $z_j =\cO(|k_j|^{-\infty})$ if for every $N>0$ there exists a $C_N>0$ such that $|z_j|\leq C_N |k_j|^{-N}$.
\end{defin} 

The significance of the assumption $\Re(k_j)\geq 1$ in Definition~\ref{def:quasi} is that it {\it(a)} bounds the wavenumbers away from zero, and {\it(b)} specifies that we are considering wavenumbers in the right-half complex plane. Therefore $\Re(k_j)\geq 1$ could be replaced by $\Re(k_j)\geq k_0$ for any $k_0>0$.

The concentration of the quasimodes near the boundary $\Gamma$ means that they are understood in the asymptotic literature as ``whispering gallery'' modes; see, e.g., \cite{BaBu:91}.

\begin{theorem}\mythmname{Existence of quasimodes \cite{PoVo:99}}\label{thm:PoVo}
If $\Oin$ is a $C^\infty$ convex domain with strictly positive curvature, $n_i>1$ and $A_N>0$, then there exists a quasimode for the transmission problem \eqref{eq:BVP2}. Furthermore, $0>\Im k_j = \cO(|k_j|^{-\infty})$ (i.e.~$k_j$ is super-algebraically close to the real axis).
\end{theorem}
The main result in \cite{PoVo:99} about resonances then follows from showing that there exists an infinite sequence of resonances that are super-algebraically close to the quasimodes \cite[Proposition~2.1]{PoVo:99}. 

Theorem \ref{thm:PoVo} implies that there exists an increasing sequence of \emph{complex} wavenumbers through which the solution operator grows super-algebraically.
We show in Corollary \ref{cor:PoVo2} below that this result implies that there exists an increasing sequence of \emph{real} wavenumbers through which one obtains this growth.

To prove this corollary we need two preparatory results. The first (Corollary \ref{cor:PoVoa} below) is that one can change the normalisation $\|u_i^{(j)}\|_{L^2(\Gamma)}=1$ in the definition of the quasimode to $\|\nabla_T u_i^{(j)}\|_{L^2(\Gamma)}+ |k_j|\|u_i^{(j)}\|_{L^2(\Gamma)}=|k_j|$; it turns out that it will be more convenient for us to work with this normalisation (note that we put the factor $|k_j|$ on the right-hand side because, since we expect $\|\nabla_T u_i^{(j)}\|_{L^2(\Gamma)}$ to be proportional to $|k_j|\|u_i^{(j)}\|_{L^2(\Gamma)}$, this normalisation therefore keeps $\|u_i^{(j)}\|_{L^2(\Gamma)}$ being $\mathcal{O}(1)$).
The second result (Lemma \ref{lem:L2}) is that, under the new normalisation, the $L^2$ norms of $u^{(j)}_{i/o}$ in $\Omega_{i/o}$ are bounded by $C|k_j|^m$ for some $C$ and $m$ independent of $j$ and $k_j$.

\begin{cor}\mythmname{Quasimode under different normalisation}\label{cor:PoVoa}
If the normalisation $\|u_i^{(j)}\|_{L^2(\Gamma)}$ $=1$ in Definition \ref{def:quasi} is changed to 
\beq\label{eq:norm}
\|\nabla_T u_i^{(j)}\|_{L^2(\Gamma)}+ |k_j|\|u_i^{(j)}\|_{L^2(\Gamma)}=|k_j|,
\eeq
then the result of Theorem \ref{thm:PoVo} still holds.
\end{cor}

This follows from the construction of the quasimode in \cite[\S5.4]{PoVo:99}; instead of dividing the $f_j$s by $\|f_j\|_{L^2(\Gamma)}$ one divides by $(\|\nabla_T f_j\|_{L^2(\Gamma)}+ |k_j|\|f_j\|_{L^2(\Gamma)})/|k_j|$. Since $\|f_j\|_{L^2(\Gamma)}$ is not super-algebraically small by \cite[last equation on page 437]{PoVo:99}, neither is $(\|\nabla_Tf_j\|_{L^2(\Gamma)}+ |k_j|\|f_j\|_{L^2(\Gamma)})/|k_j|$.

\ble\mythmname{Bound on the $L^2$ norms of $u^{(j)}_{i/o}$ in $\Omega_{i/o}$}\label{lem:L2}
The quasimode of Theorem \ref{thm:PoVo} (under the normalisation \eqref{eq:norm} in Corollary \ref{cor:PoVoa}) satisfies
\beq\label{eq:L2}
\big\|u^{(j)}_i\big\|_{L^2(\Oi)}+\big\|u^{(j)}_o\big\|_{L^2(\Oout)}\leq C  |k_j|^{1/2}+ \cO(|k_j|^{-\infty}),
\eeq
where $C$ is independent of $j$ and $k_j$.
\ele

\bpf
The plan is to obtain $k_j$-explicit bounds on the Cauchy data of $u^{(j)}_i$ and  $u^{(j)}_o$ and then use Green's integral representation and $k_j$-explicit bounds on layer potentials.
In the proof we use the notation that $a\lesssim b$ if $a\leq C b$ for some $C$ independent of $j$ and $k_j$ (but not necessarily independent of $n_i$ and $A_N$).

By the transmission condition \eqref{eq:QM3} and the normalisation \eqref{eq:norm}, we have
\begin{align*}\nonumber
\big\|\nabla_T u_o^{(j)}\big\|_{L^2(\Gamma)}+ |k_j|\big\|u_o^{(j)}\big\|_{L^2(\Gamma)}
&\leq \big\|\nabla_T  u_i^{(j)}\big\|_{L^2(\Gamma)}+ |k_j|\big\|u_i^{(j)}\big\|_{L^2(\Gamma)}+\cO(|k_j|^{-\infty})\\
&= |k_j|+\cO(|k_j|^{-\infty}).
\end{align*}

The bound in \cite[Lemma 5]{MoLu:68} on the Dirichlet-to-Neumann map for $\Oi$ that are star-shaped with respect to a ball (and thus, in particular, smooth convex $\Oi$) implies that
\begin{align}\nonumber
\big\|\dn u_o^{(j)}\big\|_{L^2(\Gamma)} &\lesssim \big\|\nabla_T u_o^{(j)}\big\|_{L^2(\Gamma)} + |k_j|\big\|u_o^{(j)}\big\|_{L^2(\Gamma)} + \cO(|k_j|^{-\infty})\\
&\lesssim |k_j| \big(1+ \cO(|k_j|^{-\infty})\big)
\label{eq:L22}
\end{align}
where the $\cO(|k_j|^{-\infty})$ in the first line is the contribution from $(\Delta +k_j^2)u_o^{(j)}$.
Note that \cite[Lemma 5]{MoLu:68} is valid for real wavenumber, but an analogous bound holds for complex wavenumbers with sufficiently small $\cO(1)$ imaginary parts -- see \cite[Theorem I.2D]{Mo:75}.
Via the transmission condition \eqref{eq:QM4}, the bound \eqref{eq:L22} holds with 
$\dn u_o^{(j)}$ replaced by $\dn u_i^{(j)}$.

The result \eqref{eq:L2} then follows from (i) Green's integral representation (applied in both $\Oout$ and $\Oin$), 
(ii) the classical bound on the free resolvent
\beq\label{eq:NP}
\N{\chi R(k) \chi}_{L^2(\Rea^d)\rightarrow L^2(\Rea^d)} \lesssim \frac{1}{|k|}\exp(a (\Im k)_-),
\eeq
for some $a>0$ (depending on $\chi$), where $x_-=0$ for $x\geq 0$ and $x_-=-x$ for $x<0$, and $\chi$ is any cutoff function, and (iii) 
the bounds on the single- and double-layer 
potentials
\beq\label{eq:SLDL}
\N{\chi \cS_k }_{L^2(\Gamma)\rightarrow L^2(\Rea^d)} \lesssim |k|^{-1/2} \exp(a (\Im k)_-)\tand
\N{\chi \cD_k }_{L^2(\Gamma)\rightarrow L^2(\Rea^d)} \lesssim |k|^{1/2}\exp(a (\Im k)_-).
\eeq
For a proof of \eqref{eq:NP} see, e.g., \cite[Theorems 3 and 4]{Va:75} or \cite[Theorem 3.1]{DyZw:16}. The bounds \eqref{eq:SLDL} with $\Im k=0$ are obtained from \eqref{eq:NP} in \cite[Lemma 4.3]{Sp2013a}; the same proof goes through with $\Im k\neq 0$.
\epf

\bre\mythmname{Sharp bounds on the single- and double-layer operators}
The bounds \eqref{eq:SLDL} are not sharp in their $k$-dependence. Sharp bounds for real $k$ are
given in \cite[Theorems 1.1, 1.3, and 1.4]{HaTa:15}, and we expect analogous bounds to be valid in the case when $k$ is complex with sufficiently-small imaginary part \cite{Ta:16}. In this case, the exponent $1/2$ in \eqref{eq:L2} would be lowered to $1/6$, but 
this would not affect the following result (Corollary \ref{cor:PoVo}), which uses \eqref{eq:L2}.
\ere

\begin{cor}\mythmname{Existence of quasimodes with real wavenumbers}\label{cor:PoVo}
Given the quasi\-mode in Theorem \ref{thm:PoVo} (under the normalisation \eqref{eq:norm}),
\beqs
\big\{\Re k_j , \big(u_i^{(j)}, u_o^{(j)}\big) \big\}_{j=1}^\infty
\eeqs
is also a quasimode (again under the normalisation in Corollary \ref{cor:PoVoa}).
\end{cor}

\bpf
From Definition \ref{def:quasi}, we only need to show that 
$\| (\Delta + n_i(\Re k_j)^2)u_i^{(j)}\|_{L^2(\Oin)}\!=\!\cO(|k_j|^{-\infty})$
and
$\| (\Delta + (\Re k_j)^2)u_o^{(j)}\|_{L^2(\Oout)}=\cO(|k_j|^{-\infty})$.
We have
\beqs
\N{ (\Delta + n_i (\Re k_j)^2)u_i^{(j)}}_{L^2(\Oi)}\leq \N{n_i\big(-2 \ri (\Re k_j) (\Im k_j) + (\Im k_j)^2\big)u_i^{(j)}}_{L^2(\Oi)} + \cO(|k_j|^{-\infty}),
\eeqs
which is $\cO(|k_j|^{-\infty})$ since $\Im k_j= \cO(|k_j|^{-\infty})$ and the $L^2$ norm of $u_i^{(j)}$ is bounded by $C|k_j|^{1/2}$ by the bound \eqref{eq:L2}. Similarly, we have $\| (\Delta + (\Re k_j)^2)u_o^{(j)}\|_{L^2(\Oout)}=\cO(|k_j|^{-\infty})$.
\epf

The super-algebraic growth of the solution operator through real values of $k$ can therefore be summarised by the following.

\begin{cor}\mythmname{Super-algebraic growth through real wavenumbers}\label{cor:PoVo2}
If $\Oin$ is a $C^\infty$ convex domain with strictly positive curvature, $n_i>1$ and $A_N>0$, then there 
exists a sequence
\beqs
\big\{ \widetilde{k}_j , \big(u_i^{(j)}, u_o^{(j)}\big) \big\}_{j=1}^\infty,
\eeqs
where $\widetilde{k}_j \in \Rea$, $\widetilde{k}_j\tendi$, $u_{i/o}^{(j)} \in C^\infty(\overline{\Omega_{i/o}})$, the support of $u^{(j)}_{i/o}$ is contained in a fixed compact neighbourhood of $\Gamma$, and for all $N\in\IN$ there exists $C_N>0$ independent of $\widetilde k_j$ such that
\begin{align*}
&\min\Big\{\big\|u_i^{(j)}\big\|_{H^{3/2}(\Oin)},\big\|u_o^{(j)}\big\|_{H^{3/2}(\Oout)}\Big\}\\
&\ge C_N \widetilde k_j^N \Big(\big\|(\Delta + n_i \widetilde{k}_j^2)u_i^{(j)}\big\|_{L^2(\Oin)}
+\big\|(\Delta +\widetilde{k}_j^2)u_o^{(j)}\big\|_{L^2(\Oout)}
\\&\hspace{20mm}
+\N{ u_i^{(j)}-u_o^{(j)} }_{H^2(\Gamma)} +\N{ \dn u_i^{(j)}-A_N \dn u_o^{(j)} }_{H^2(\Gamma)}
\Big).
\end{align*}
\end{cor}
\bpf
This follows from Corollary \ref{cor:PoVo} by the trace theorem (see, e.g., \cite[Theorem 3.37]{MCL00}).
\epf

\subsection{The results of \texorpdfstring{\cite{Cap12,CLP12}}{Capdeboscq} when \texorpdfstring{$\Oi$}{Omega-i} is a ball and numerical examples}
\label{sec:ball}

The results of \cite{Cap12, CLP12} (summarised in \cite[Chapter 5]{AC16}) consider the case when $\Oi$ is a 2- or 3-d ball. It is convenient to discuss these alongside numerical examples in 2-d. (We therefore only discuss the 2-d results of \cite{Cap12}, but the 3-d analogues of these are in \cite{CLP12}.)

When $\Oi$ is the 2-d unit ball ($\Oin=B_1$), $u_i$ and $u_o$ can be expressed in terms of Fourier series $\ee^{\ri \nu \theta}$, where $\nu\in \mathbb{Z}$ and $\theta$ is the angular polar coordinate, with coefficients given in terms of Bessel and Hankel functions; see, e.g., \cite[\S12]{Ga:15} or \cite[\S5.3]{AC16}. 
The resonances of the BVP \eqref{eq:BVP2} (with $n_i>1$) are then the (complex) zeros of 
\beqs
F_{\nu}(k):=\AN\sqrt{n_i} \,J'_\nu(\sqrt{n_i}\, k)H^{(1)}_\nu (k) -  H^{(1)\prime}_\nu (k) J_\nu(\sqrt{n_i}\,k),
\quad \nu\in\IZ,
\eeqs
see, e.g., \cite[Equation (81)]{Ga:15}. 
The asymptotics of the zeros of $F_{\nu}(k)$ in terms of $\nu$ are given by, e.g., \cite{LaLeYo:92, Sc:93}.
Indeed, let $k_{\nu,m}$ denote the $m$th zero of $F_{\nu}(k)$ with positive real part (where the zeros are ordered in terms of magnitude), and let $\alpha_m$ denote the $m$th zero of the Airy function ${\rm Ai}(-z)$. 
Then by \cite[Equation~(1.1)]{LaLeYo:92}, for fixed $m$,
\beq\label{eq:resonance_asym}\sqrt{n_i}\, k_{\nu,m} = \nu + 2^{-1/3}\alpha_m \nu^{1/3} + \cO(1)\quad\tas\quad \nu\tendi;
\eeq
recall that the result \cite[Proposition~2.1]{PoVo:99} (discussed in \S\ref{sec:PoVo}) implies that the resonances are real to all algebraic orders.

\paragraph{Exponential growth in $k$.}
For our first numerical examples, we take $\nin>1$, $\nout=1$, $\AN=1$ and for $\nu\in\IN_0$ we let $k$ be the real part of the first resonance corresponding to the angular dependence $\ee^{\ri \nu\theta}$ (i.e.~$k=\Re k_{\nu,1}$ in the previous paragraph).
We then let $f_o=0$ and $f_i= c_\nu J_{\nu}(kr)\ee^{\ri \nu \theta}$, where $c_\nu=(\pi(J_\nu^2(k)-J_{\nu+1}(k)J_{\nu-1}(k)))^{-1/2}$ 
is such that $\|f_i\|_{L^2(\Oin)}=1$ for all $\nu\in \IN_0$.
The field $u$ is then
\begin{align*}
u_i&=c_\nu\frac{\ee^{\ri\nu\theta}}{k^2(n_i\!-\!1)}\bigg(
J_\nu(kr)-J_\nu(k\sqrt{\nin}r)
\frac{\AN J_\nu'(k)\Hankel_\nu(k)-H^{(1)'}_\nu(k) J_\nu(k)}
{F_\nu(k)}
\bigg),
\\
u_o&=c_\nu\frac{\ee^{\ri\nu\theta}\Hankel_\nu(kr)}{k^2(n_i\!-\!1) \Hankel_\nu(k)}
\bigg(J_\nu(k)-J_\nu(k\sqrt{\nin})
\frac{\AN J_\nu'(k)\Hankel_\nu(k)-H^{(1)'}_\nu(k) J_\nu(k)}
{F_\nu(k)}
\bigg).
\end{align*}
Figure \ref{fig:1} plots the $k\sqrt n$-weighted $L^2$ norm and the $H^1$ seminorm of the particular $u_{i/o}$ above in $\Oi$ and $D_R$ respectively (with $R$ taken to be 2) as $\nu$ runs from $0$ to $\nu_{\max}$. 
The left and right panels in Figure~\ref{fig:1} show the norms of $u_{i/o}$ for $\nin=3$ and $\nu_{\max}=64$, and for $\nin=10$ and $\nu_{\max}=23$, respectively.
We see that the norms of the solution appear to grow exponentially; the existence of a sequence growing super-algebraically is expected by Corollary \ref{cor:PoVo2}, although we are considering lower-order norms than in this result.
We note that the estimate \eqref{eq:resonance_asym} is not enough to deduce exponential growth of the scattered field for $k=\Re k_{\nu,1}$ and increasing $\nu$, but a more refined analysis is needed (taking into account the fact that the imaginary parts of the resonances are superalgebraically small).
The result \cite[Theorem 5.4]{AC16} proves that, at least in the case of plane-wave incidence, Sobolev norms of $u_{o}$ on spherical surfaces in $D_R$ sufficiently close to $\Gamma$ grow exponentially through a sequence of real wavenumbers, where the wavenumbers are defined in terms of the resonances by \cite[Lemma 5.16]{AC16}; this result can then be used to prove exponential growth in the $L^2$ norm of $u_o$.

\begin{figure}[htb]
\centering    
\includegraphics[width=.49\textwidth, clip, trim={50 0 45 0}]{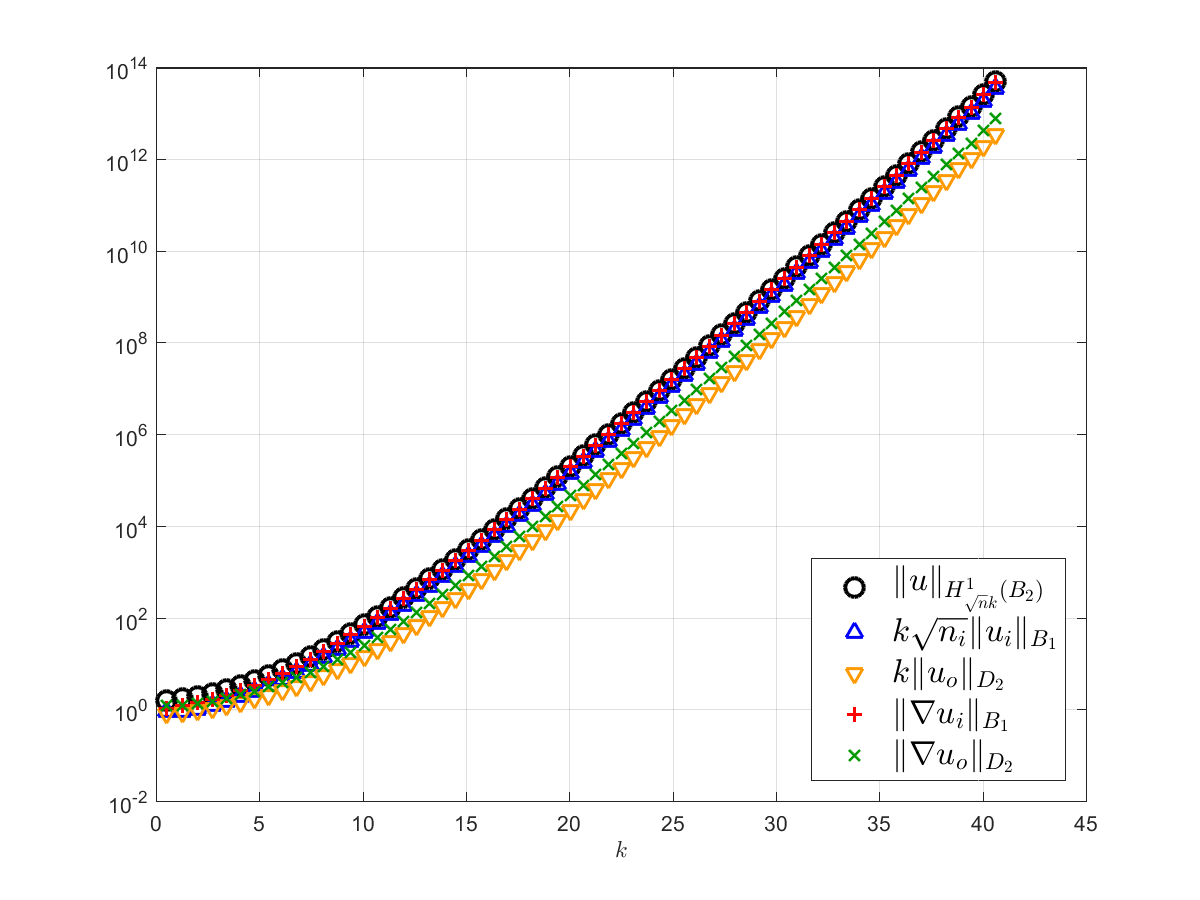}
\includegraphics[width=.49\textwidth, clip, trim={50 0 45 0}]{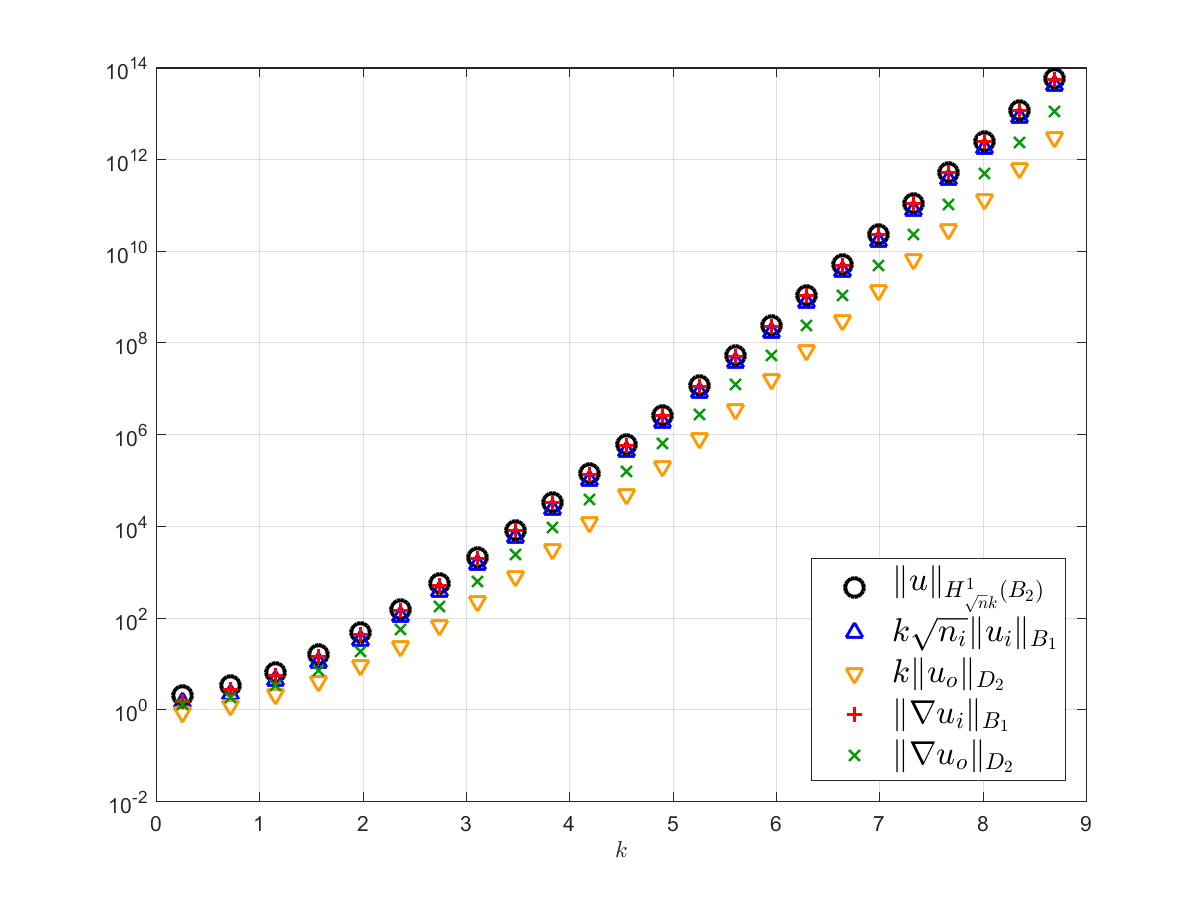}
\caption{Left plot: the norms of the solution of BVP \eqref{eq:BVP2} with $\Oin=B_1$, $\nin=3,\AN=1$, $\N{f_i}_\Oin=1$, $f_o=0$, and where each $k$ (on the abscissa) is the real part of the first resonance whose angular dependence is $\ee^{\ri \nu\theta}$ for $\nu=0,\ldots,64$, i.e.\ $k=\Re k_{\nu,1}$.
The black circles show the weighted $H^1(B_2)$ norms of $u$ as in the left-hand side of \eqref{eq:bound1}, while the other four markers represent the weighted $L^2$ norms and the $H^1$ seminorms of $u_i$ and $u_o$ in $\Oi=B_1$ and $D_2=B_2\setminus\overline{B_1}$.
\newline Right plot: same as above with $\nin=10$ and $\nu=0,\ldots,23$.
}
\label{fig:1}
\end{figure}

\begin{figure}[htb]
\centering
\includegraphics[width=70mm, clip, trim={80 30 45 30}]
{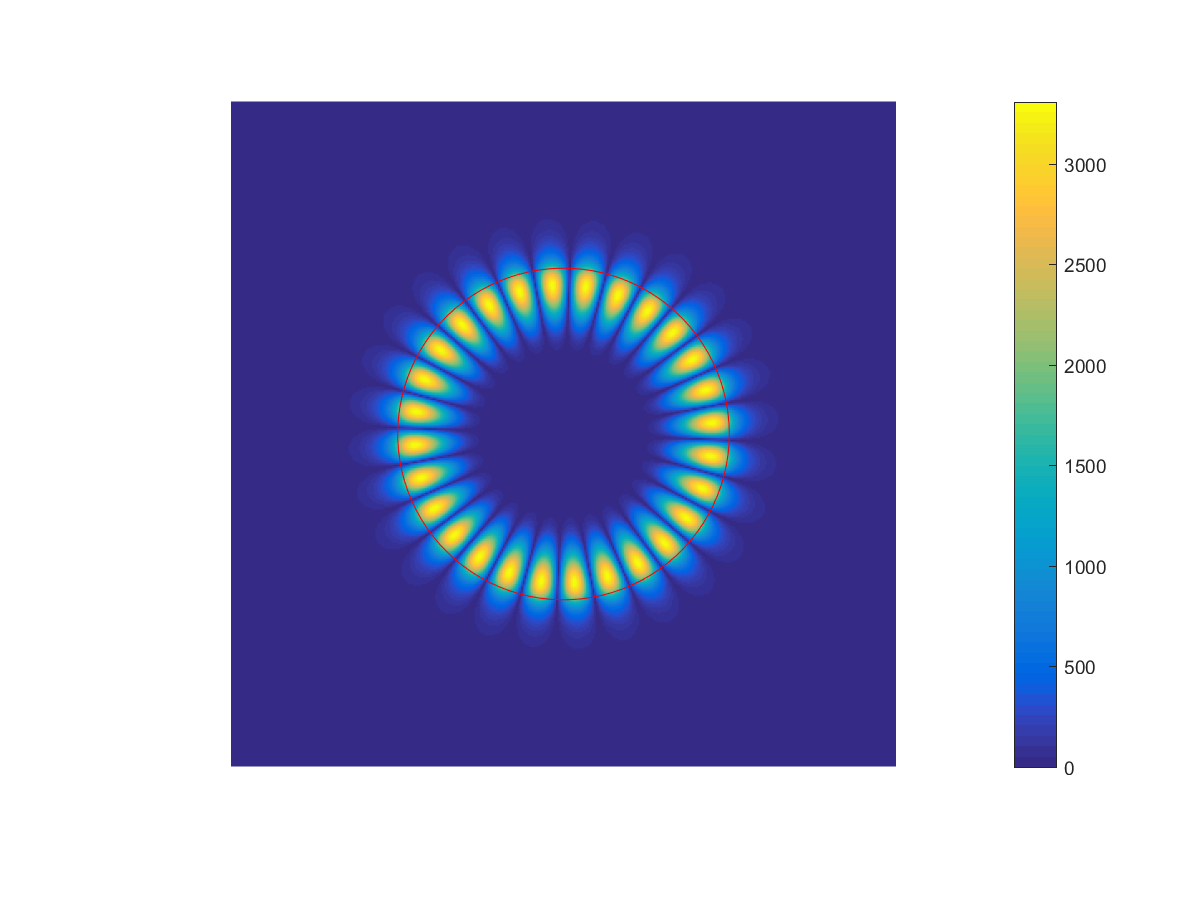}
\includegraphics[width=70mm, clip, trim={80 30 45 30}]
{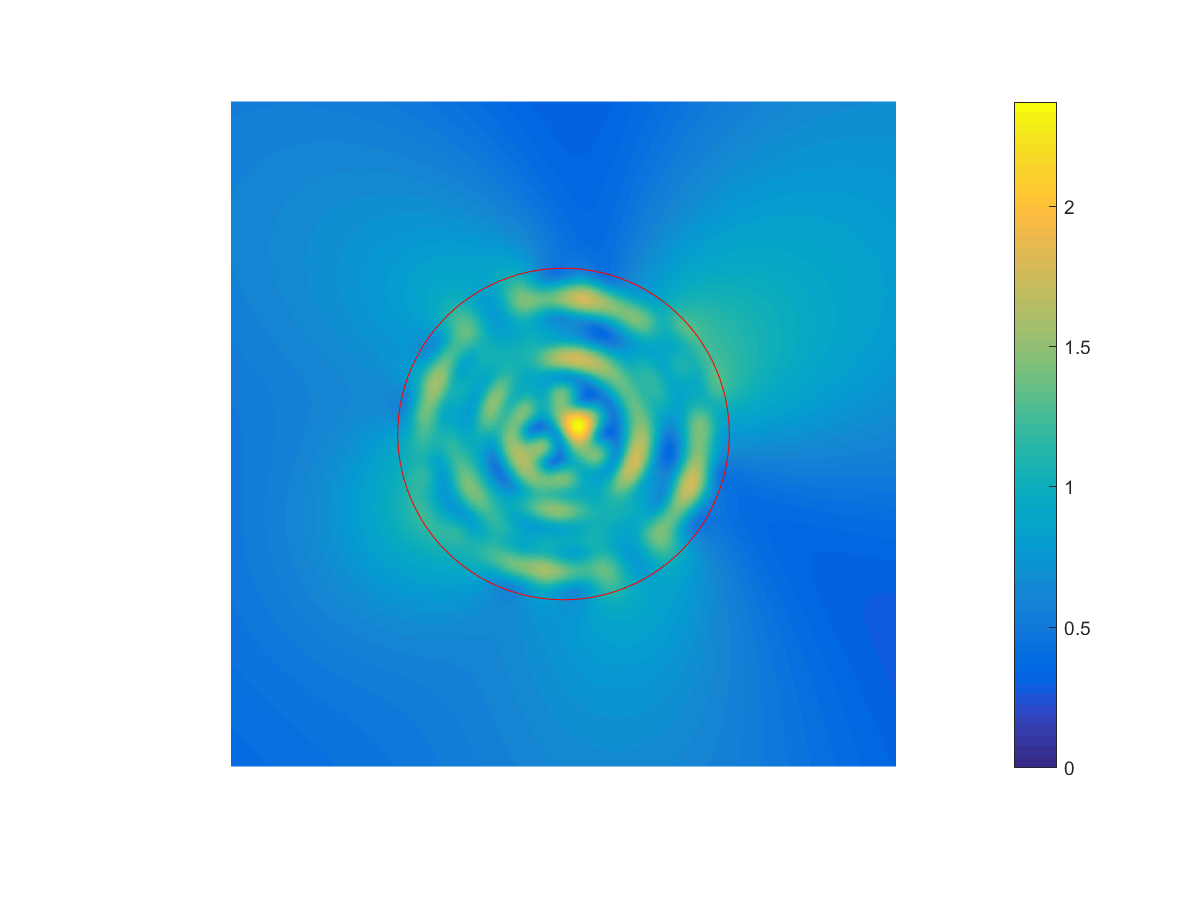}
$k_1=1.77945199481921\approx\Re k_{14,1}$\hspace{28mm} $k_3=1.779451994815$\hspace{20mm}

\includegraphics[width=70mm, clip, trim={80 30 45 30}]
{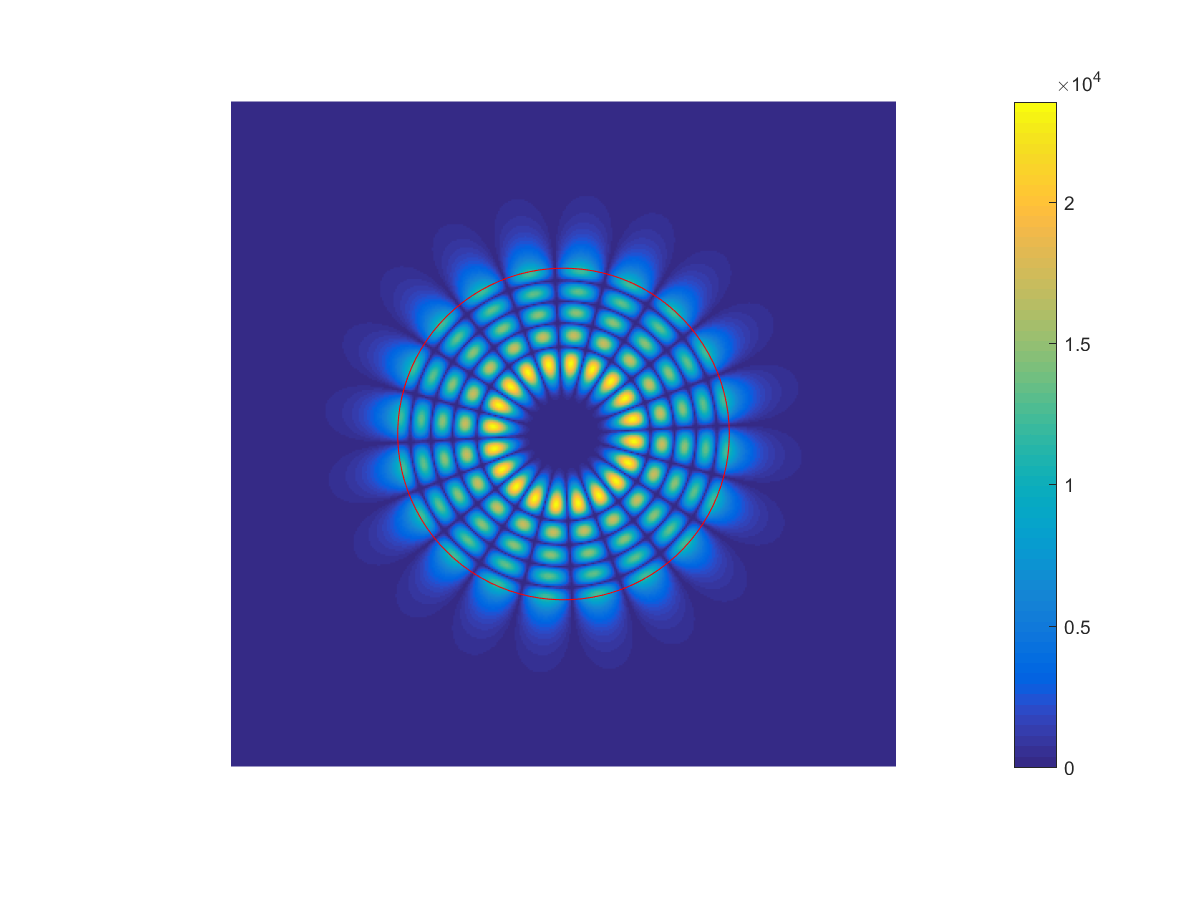}
\includegraphics[width=70mm, clip, trim={80 30 45 30}]
{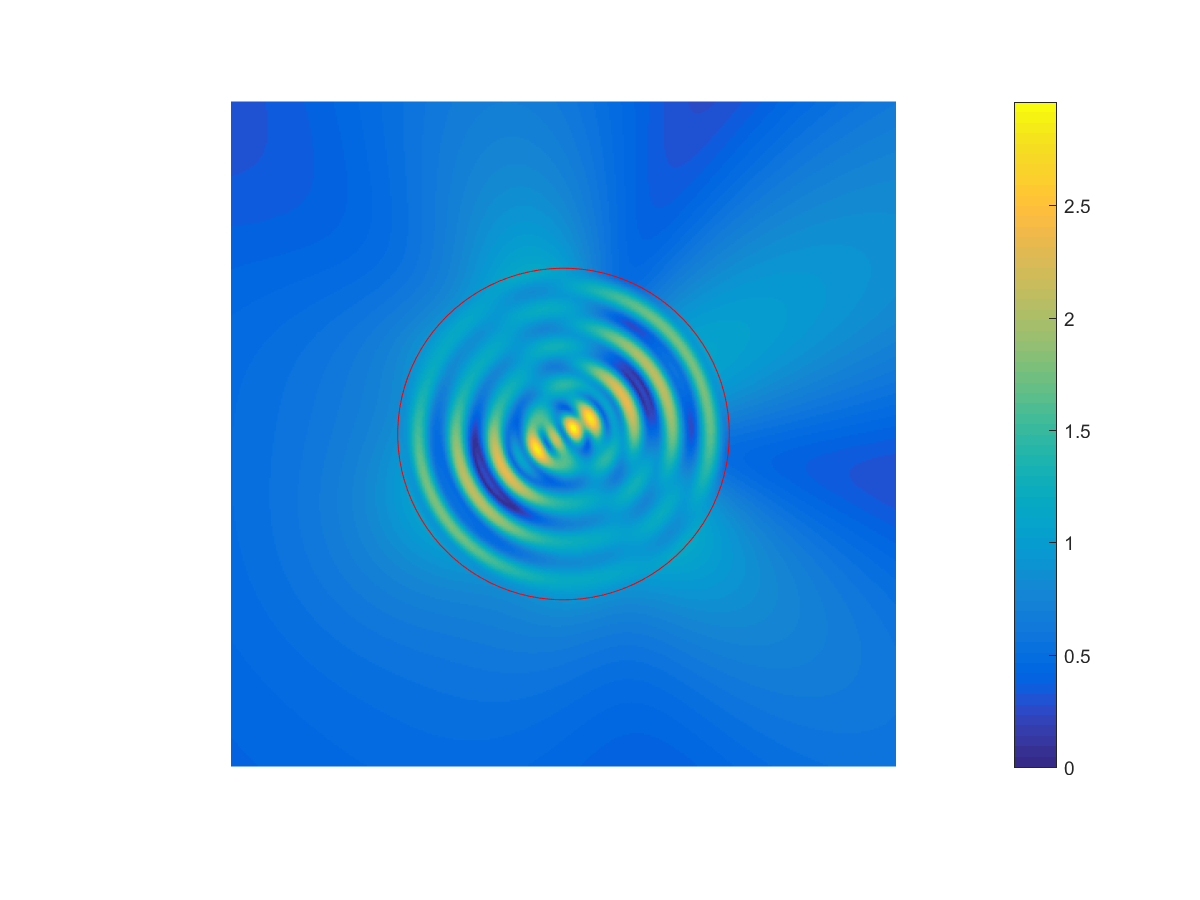}
$k_2=2.75679178324354\approx\Re k_{10,5}$\hspace{35mm} $k_4=2.757$\hspace{25mm}
\caption{The absolute value of the field scattered by plane waves $\ee^{\ri k(x\cos\phi,y\sin\phi)}$ with unit amplitude and propagation direction $\phi=\pi/6$ hitting the disc $B_1$ with $\nin=100$.
The red line denotes the boundary $\Gamma$ of $\Oin=B_1$.
The four figures differ only in the wavenumber.
\newline
In the upper left plot, for $k_1\approx\Re k_{14,1}$ the impinging plane wave excites a typical whispering gallery mode, quickly decaying away from the interface $\Gamma$.
In the lower left plot, $k_2\approx\Re k_{10,5}$, the absolute value of the excited mode has $m=5$ peaks in the radial direction and $2\nu=20$ in the angular direction.
In the upper right plot $k_3$ differs from $k_1$ only by a factor of order $10^{-12}$ but nevertheless generates a completely different plot, as the quasimode is not excited; note also the different scales shown in the colour bar.
In the lower right plot, $k_4$ differs from $k_2$ only by a factor of order $10^{-4}$, and the same considerations apply.
}
\label{fig:Fields}
\end{figure}

\paragraph{Localisation of $u_o$ in $\Oout$ at resonant frequencies.}

The plots in Figure \ref{fig:Fields} show the absolute value of the fields scattered by a plane wave impinging on the unit disc with $\nin=100$.
In the left plots, a wavenumber equal to the real part of a resonance $k_{\nu,m}$ excites a quasimode. 
In both these examples, $u_{o}$ is localised close to $\Oin$; this is expected both from Theorem~\ref{thm:PoVo} above, and from 
\cite[Theorem 5.2]{AC16}. Indeed, this latter result gives bounds on $u_{o}$ for all values of $k$ and $\nin$, but if $\nin>1$ they hold only in the ``far field'', i.e.\ at distance at least $(\nin-1)\diam\Oin/2$ from $\Oin$, showing that the quasimodes generate large fields in a small neighbourhood of $\Oin$ only.

\paragraph{Sensitivity to the wavenumber.}

The right plots in Figure \ref{fig:Fields} show how a small perturbation of the wavenumber $k$ (e.g.\ by a relative factor of about $2.4\cdot10^{-12}$ from that in the upper left plot to that in the upper right one, and of $7.6\cdot10^{-5}$ in the second row) avoids the quasimode and gives a solution $u$ with much smaller norm; see the scale displayed by the colour bar.

This phenomenon suggests that for certain data the exponential blow-up of the solution operator 
when $\nin> 1$ can be avoided. 
Indeed, given $k$, if $f_{i/o}(r,\theta)=\sum_\ell \alpha_{\nu,i/o}(r)\ee^{\ri \nu \theta}$ and $\alpha_{\nu,i/o}\equiv 0$ for any $\nu$ such that there exists an $m$ such that $k_{\nu,m}\approx k$, then 
the resonant modes should be excluded from the solution.
In the context of scattering by an incident field, \cite[Theorem~6.5]{Cap12} describes the size of the neighbourhood $I$ of $\{\Re k_{\nu,m}\}_{\nu\in\IZ,m\in\IN}$, such that the scattered field is uniformly bounded for $k\in(0,\infty)\setminus I$ (see also \cite[Remark 5.5]{AC16}).
The strong sensitivity of the quasimodes with respect to the wavenumber explains why their presence can go unnoticed even by extensive numerical calculation; see, e.g.,\cite{BaChGo:16} for a related class of problems.

\appendix

\section{Proof of Lemma \texorpdfstring{\ref{lem:exist}}{2.7} (existence, uniqueness, and regularity of the BVP solution)}\label{sec:appA}

\paragraph{Uniqueness:} 
In this setting of Lipschitz $\Oi$, uniqueness follows from Green's identity and a classical result of Rellich (given in, e.g., \cite[Lemma 3.11 and Theorem 3.12]{CoKr:83}); the proof for our assumptions on the parameters $k, n_i,n_o,a_o,a_i, A_D,$ and $A_N$ is short, and so we give it here.

With $f_i=0$, $f_o=0$, $g_D=0$, $g_N=0$, we apply Green's identity to $u_i$ in $\Oi$ (this is allowed since $\Oi$ is Lipschitz and $u_i \in H^1(\Oi,\Delta)$ by, e.g., \cite[Theorem 4.4]{MCL00}) to obtain
\beqs
\int_\Gamma a_i u_i \overline{\dn u_i}  - a_i \int_{\Oin}\ngus + n_i \overline{k^2} \int_\Oin \nus=0.
\eeqs
Using the transmission conditions in \eqref{eq:BVP} we get 
\beq\label{eq:A1}
\frac{a_o}{A_D A_N}\int_\Gamma u_o\overline{\dn u_o}  - a_i \int_{\Oin}\ngus + n_i \overline{k^2} \int_\Oin \nus=0.
\eeq
By \cite[Theorem 3.12]{CoKr:83}, if 
\beqs
\Im\left( k \sqrt{\frac{n_o}{a_o}} \int_\Gamma u_o \overline{\dn u_o}\right)\geq 0
\eeqs
then $u_o\equiv 0$ in $\Oout$; the Cauchy data of $u_i$ is then zero by the transmission conditions, and thus $u_i\equiv 0$ in $\Oin$. 
Note that \cite[Theorem 3.12]{CoKr:83} is stated and proved for $d=3$, but the proof goes through in an identical way for $d=2$ and $d\geq 4$ using the radiation condition \eqref{eq:src} and the asymptotics of the fundamental solution in these dimensions (see, e.g., \cite[Theorem~3.5]{DyZw:16}).

Multiplying \eqref{eq:A1} by $k\sqrt{n_o/a_o} (A_N A_D/a_o)$ and recalling that all the parameters apart from $k$ are real, we see that a sufficient condition for uniqueness is 
\beqs
\Im\left( k a_i\int_{\Oin}\ngus - k n_i \overline{k^2} \int_\Oin \nus\right)\geq 0;
\eeqs
this inequality holds since $n_i$ and $a_i$ are real, and $\Im k\geq 0$. 

\paragraph{Existence and regularity:} 

At least in the case $g_D=g_N=0$, existence immediately follows from uniqueness via Fredholm theory. Here, however, we use the integral-equation argument of \cite{ToWe:93}, since this also establishes the regularity results on $\Gamma$. This integral-equation argument is the Lipschitz analogue of the argument for 
 sufficiently smooth $\Gamma$ in \cite[Theorem 4.6]{KrRo:78} (see also \cite[Corollary 4.6]{CoSt:85}, which  covers 2-d polygons).

If $f_i=0$ and $f_o=0$ then the existence and the regularity of $u$ follow from the integral-equation argument in \cite[Theorem 7.2]{ToWe:93};
 to match their notation we choose $u_e=\frac{\aout}{\AN\ain}u_o$ and
\begin{align*}
k_1=k\sqrt{\frac\nin\ain},\quad
k_2=k\sqrt{\frac\nout\aout},\quad
\mu_1=\AD\aout,\quad
\mu_2=\AN\ain,\quad
f=\aout g_D,\quad
g=\frac1{\AN\ain}g_N.
\end{align*}
The result \cite[Theorem 7.2]{ToWe:93} is stated only for $d= 3$ and $k\in \Rea\setminus\{0\}$, but we now outline why it also holds for $d\geq 2$ and $\Im k>0$; we first discuss the dimension. 

The two reasons Torres and Welland only consider $d=3$ is that 
\ben
\item Their argument treats the Helmholtz integral operators as perturbations of the corresponding Laplace ones, and there 
is a technical difficulty that the fundamental solution of Laplace's equation does not tend to zero at infinity when $d=2$, whereas it does for $d\geq 3$.
\item The case $d \geq 4$ is very similar to $d=3$, except that the bounds on the kernels of the integral operators are slightly more involved.
\een

Regarding the second reason: the only places in \cite{ToWe:93} where these bounds are used and that contribute to the existence result \cite[Theorem 7.2]{ToWe:93} are Part (vi) of Lemma 6.2, and Points (i)--(vi) of Page 1466. The analogue of the Hankel-function bounds in \cite[Equation 2.25]{CGLS12} for $d\geq 4$ can be used to show that these results in \cite{ToWe:93} hold for $d\geq 4$. 

Regarding the first reason: the harmonic-analysis results about Laplace integral operators, upon which Torres and Welland's proof rests, also hold when $d=2$. More specifically, the results in Sections 4 and 6 of \cite{ToWe:93} hold when $d=2$ by \cite{Ve:84} and \cite{CoMcMe:82} (a convenient summary of these harmonic-analysis results is given in \cite[Chapter 2]{CGLS12}). The results (i) and (iii)--(vi) on Page 1463 of \cite{ToWe:93} hold when $d=2$ by results in \cite[\S4]{Ve:84}, and the result (ii) holds if one defines the Laplace fundamental solution as $(1/2\pi)\log(a/|\bx-\by|)$ where the constant $a$ is not equal to the so-called ``capacity'' of $\Gamma$; see \cite[Page 115]{CGLS12}, \cite[Theorem 8.16]{MCL00}. 
Finally, Lemma 3.1 of \cite{ToWe:93} holds when $d=2$ by \cite[Theorem~2]{EsFaVe:92}.
All this means that Lemma 3.2 of \cite{ToWe:93} holds when $d=2$, which in turn means that \cite[Theorem 7.2]{ToWe:93} holds when $d=2$.

We therefore have that \cite[Theorem 7.2]{ToWe:93} holds for $d\geq 2$ and for $k\in \Rea\setminus\{0\}$. Inspecting the proof of \cite[Theorem 7.2]{ToWe:93} we see that the properties of the Helmholtz boundary-integral operators used are unchanged if $\Im k>0$. Indeed, the compactness of the differences of the Helmholtz and Laplace boundary-integral operators holds for $\Im k>0$ by (v) and (vi) on Page 1466 of \cite{ToWe:93} and by (vi) in Lemma 6.2 of \cite{ToWe:93}, and the uniqueness of the BVP (used on Page 1483) holds by the argument above. Therefore \cite[Theorem 7.2]{ToWe:93} holds for $d\geq 2$ and for $k\in\Com\setminus\{0\}$ with $\Im k\geq 0$.

If the volume source terms $f_i$ and $f_o$ are different from zero, we define $w_i$ and $w_o$ to be the solutions to the following problems:
\begin{align*}
\begin{cases}
w_i\in H^1(\Oin)\\
\ain\Delta w_i+k^2\nin w_i=f_i &\iin \Oin,\\
\dn w_i-\ri k w_i=0 &\oon\Gamma.
\end{cases}
\qquad
\begin{cases}
w_o\in H^1\loc(\Oout)\\
\aout\Delta w_o+k^2\nout w_o=f_o &\iin \Oout,\\
w_o=0 &\oon\Gamma,\\
w_o \in \SRC(k\sqrt{n_o/a_o}).
\end{cases}
\end{align*}
Then by \cite[Proposition~3.2]{MOS12} $w_i\in H^1(\Oin)$, $\dn w_i\in L^2(\Gamma)$ and $w_i\in H^1(\Gamma)$, and by \cite[Part (i) of Lemma 3.5]{Sp2013a}
$w_o\in H^1\loc(\Oout)$, $\dn w_o\in L^2(\Gamma)$ and $w_o\in H^1(\Gamma)$
(note that both of the results \cite[Proposition~3.2]{MOS12} and \cite[Part (i) of Lemma 3.5]{Sp2013a} are specialisations of the regularity results of Ne\v{c}as \cite[\S5.1.2 and \S5.2.1]{Ne:67} to Helmholtz BVPs).
Then $\widetilde u_i:=u_i-w_i$ and $\widetilde u_o:=\frac{\aout}{\AN\ain}(u_o-w_o)$ satisfy problem (P) of \cite{ToWe:93} with $k_1,k_2,\mu_1,\mu_2$ as above and 
$$f=\aout g_D+\AD\aout w_i\in H^1(\Gamma), \qquad g=\frac{g_N-\aout\dn w_o}{\AN\ain}+\dn w_i\in L^2(\Gamma).$$
The existence and regularity of $u_i,u_o$ follow by applying again Theorem 7.2 of \cite{ToWe:93} to  $\widetilde u_i,\widetilde u_o$.

\section*{Acknowledgements}
Over the last few years, several colleagues have asked the authors questions about wave\-number-explicit bounds on Helmholtz transmission problems; these include Leslie Greengard (New York University), Samuel Groth (MIT), David Hewett (University College London), and Ralf Hiptmair (ETH Z\"urich).
We thank Giovanni S.\ Alberti (Genova), Simon Chandler-Wilde (University of Reading), Jeffrey Galkowski (McGill University), Ivan Graham (University of Bath), Zo\"is Motier (Rennes), Melissa Tacy (Australian National University), and Georgi Vodev (Universit\'e de Nantes) for useful discussions. 
A.~Moiola acknowledges support from EPSRC through the project EP/N019407/1, from GNCS/INDAM, and from MIUR through the ``Dipartimenti di Eccellenza'' Programme (2018--2022)--Dept. of Mathematics, Pavia.

\let\oldbibliography\thebibliography
\renewcommand{\thebibliography}[1]{\oldbibliography{#1}\setlength{\itemsep}{0pt}}

\addcontentsline{toc}{section}{References}
%\bibliographystyle{siam}
%\bibliography{references_andrea2}

\end{document}